\documentclass[smallextended,envcountsect]{svjour3}
\smartqed
\usepackage{graphicx}

\usepackage{amsfonts,amssymb,amsmath}
\journalname{JOTA}

\DeclareMathOperator*{\argmin}{argmin}

\begin{document}

\title{Minimax and Viscosity Solutions of Hamilton-Jacobi-Bellman Equations for Time-Delay Systems}


\author{Anton Plaksin}

\institute{Anton Plaksin \at
              N.N. Krasovskii Institute of Mathematics and Mechanics (IMM UB RAS)\\
              Ural Federal University\\
              Yekaterinburg, Russia,\\
              a.r.plaksin@gmail.com
}

\date{Received: date / Accepted: date}

\maketitle

\begin{abstract}
The paper deals with a Bolza optimal control problem for a dynamical system which motion is described by a delay differential equation under an initial condition defined by a piecewise continuous function. For the value functional in this problem, the Cauchy problem for the Hamilton-Jacobi-Bellman equation with coinvariant derivatives is considered. Minimax and viscosity solutions of this problem are studied. It is proved that both of these solutions exist, are unique and coincide with the value functional.
\end{abstract}


\section{Introduction}

In optimal control problems for dynamical systems which motions are described by ordinary differential equations, studies of infinitesimal properties of a value function lead to a Hamilton-Jacobi-Bellman (HJB) equation, which is a particular case of Hamilton-Jacobi (HJ) equations with partial derivatives. In the case when an optimal control problem is considered on a finite time interval and has a cost functional of Bolza type, a value function satisfies the corresponding natural terminal condition, which, together with the HJB equation, determine the Cauchy problem. Since in many cases Cauchy problems for HJ equations do not have a classical (continuously differentiable) solution, various approaches to a notion of a generalized solution were developed. The main of them are minimax and viscosity approaches. The minimax approach \cite{Subbotin_1980,Subbotin_1984,Subbotin_1995} originates in the positional differential game theory \cite{Krasovskii_Subbotin_1988,Krasovskii_Krasovskii_1995}. According to this approach, a generalized (minimax) solution is a function that satisfies the pair of stability conditions with respect to characteristic differential inclusions. In infinitesimal form, these conditions reduce to the pair of inequalities for directional derivatives. In the viscosity approach \cite{Crandall_Lions_1983,Crandall_Evans_Lions_1984}, a HJ equation is replaced by the pair of inequalities for sub- and supergradients, and a generalized (viscosity) solution is a function satisfying these inequalities. In investigations of minimax and viscosity solutions of Cauchy problems for HJB equations, it was shown (see, e.g., \cite{Subbotin_1995,Barbu_1986,Bardi_Capuzzo-Dolcetta_1997,Evans_1998}) that both of these solutions exist, are unique and coincide with the value function in the corresponding optimal control problems. The goal of the paper is to obtain the similar result in the case when a motion of a dynamical system is described by delay differential equations.

The first investigations of control problems for time-delay systems showed (see \cite{Krasovskii_1962,Osipov_1971} and also \cite{Oguztoreli_1966,Banks_1968,Banks_Manitius_1974}) that an analogue of a value function in such problems is a value functional on a space of motion histories.
It raises natural questions about the suitable notion of the differentiability of such functionals, the corresponding notions of directional derivatives, sub- and supergradients, and  definitions of generalized solutions of the corresponding HJB equations.

The viscosity solution theory for HJ equations with Frechet derivatives began with \cite{Crandall_Lions_1985,Crandall_Lions_1986a}. In these papers, the definition of the viscosity solution in terms of inequalities for Frechet sub- and supergradients was given, and existence and uniqueness of such solution were proved. After that, a lot of investigations (see, e.g., \cite{Barbu_Barron_Jensen_1988,Soner_1988,Cannarsa_Da_Prato_1990,Cannarsa_Frankowska_1992,Li_1995}) dealt with applications of the viscosity approach to control problems for abstract evolution systems in Hilbert or Banach spaces. In particular, in \cite{Soner_1988,Cannarsa_Da_Prato_1990}, for Bolza optimal control problems for evolution systems, modified definitions of viscosity solutions of Cauchy problem for HJB equations were given, their existence, uniqueness and coincidence with the value functional were shown. Note that the conditions in these papers allow to interpret some class of time-delay systems as evolution systems, however, this class is not general enough, since it does not contain systems with discrete delay. The optimal control problem for systems with discrete delay was considered in \cite{Barron_1990}. It was proved that the value functional is a viscosity solution of the Cauchy problem for the HJB equation, but the uniqueness question of the viscosity solution was not investigated. One could also mention papers \cite{Wolenski_1994,Clarke_Wolenski_1996} in which optimization problems for quite general delay differential inclusions (which cover the case of discrete delay) were considered and various necessary optimality conditions were given.

In \cite{Kim_1999}, for the description of infinitesimal properties of a value functional in optimal control problems for time-delay systems, the notion of coinvariant derivatives was used. Note that such derivatives and their close analogues were applied later to a wide range of control problems for various functional differential systems (see, e.g., \cite{Lukoyanov_2000,Lukoyanov_2001,Lukoyanov_2010a,Lukoyanov_2010b,Aubin_Haddad_2002,Pepe_Ito_2012,Lukoyanov_Gomoyunov_Plaksin_2017,Bayraktar_Keller_2018}). The theory of minimax and viscosity solutions of Cauchy problems for HJ equations with coinvariant derivatives and its application to differential games for time-delay systems were developed in \cite{Lukoyanov_2000,Lukoyanov_2001,Lukoyanov_2010a,Lukoyanov_2010b}. In these papers, the class of time-delay systems under consideration is quite general and includes systems with discrete delay. In \cite{Lukoyanov_2000,Lukoyanov_2001}, it was shown that the value functional is the unique minimax solution. In \cite{Lukoyanov_2010a}, the description of the value functional in terms of suitable directional derivatives was given. In \cite{Lukoyanov_2010b}, similar to \cite{Soner_1988}, the modified definition of a viscosity solution based on a sequence of compact sets is considered. It allows to prove that the viscosity solution exists, is unique and coincides with the minimax solution, however, such definition is not reduced to the classical definition of a viscosity solution in the particular case without delay. For more natural definitions of a viscosity solution, the uniqueness questions is still open.

This paper is aimed to solve this question and to develop the theory of minimax and viscosity solutions of HJB equations for time-delay systems, which generalizes in a natural way the classical theory of both minimax and viscosity solutions of HJB equations for systems of ordinary differential equations.

In the paper, a Bolza optimal control problem for a time-delay system with discrete delay is considered. For the value functional of this problem, a HJB equation with coinvariant derivatives is investigated. Definitions of minimax and viscosity solutions (which are consistent with the classical definitions) of the Cauchy problems for this equation are studied. It is proved that both of these solutions exist, are unique and coincide with the value functional. Besides, the feedback scheme for constructing the optimal control by the minimax (viscosity) solution is given (see the proof of Theorem \ref{thm2} $(b) \Rightarrow (a)$).

A principle idea for obtaining these results is to use the space of piecewise continuous functions as the state space in which the optimal control problem and the HJB equation are considered. As already noted earlier (see \cite{Wolenski_1994}), the choice of a suitable state space plays an important role for an application of the viscosity approach to HJB equations for time-delay systems. The space of measurable functions can be used as the state space. But such choice significantly narrows the class of the corresponding time-delay systems and excludes important for applications systems with discrete delay (see \cite{Soner_1988}). The space of continuous functions can also be used as the state space. It allows to cover the case of systems with discrete delay, but, as mentioned above, it makes it possible to prove the uniqueness only of the modified viscosity solution (see \cite{Lukoyanov_2010b}). In \cite{Barron_1990,Wolenski_1994}, other functional spaces were considered as the state space, but the uniqueness question of the viscosity solutions was not investigated.
Presented in this paper choice of the space of piecewise continuous functions allows on the one hand to consider the case of systems with discrete delays, and on the other hand, to prove the uniqueness of the viscosity solutions in the classical sense. Note that this proof is based on Lemma~\ref{lem_MVI}, which is an analogue of the theorem about "Mean value inequality" \cite{Clarke_Ledyaev_1994,Clarke_Ledyaev_Stern_Wolenski_1998} (see also \cite{Subbotin_1993}) for functionals defined on the space of piecewise continuous functions.

\vspace*{-0.2cm}

\section{Formulation of Results}

Let $\mathbb R^n$ be the $n$-dimensional Euclidian space with the inner product $\langle \cdot, \cdot \rangle$ and the norm $\|\cdot\|$. A function $x(\cdot) \colon [a,b] \mapsto \mathbb R^n$ is called piecewise continuous if there exist numbers $a = \xi_1 < \xi_2 < \ldots < \xi_k = b$ such that, for each $i \in \overline{1,k-1}$, the function $x(\cdot)$ is continuous on the interval $[\xi_i,\xi_{i+1})$ and there exist a finite limit of $x(\xi)$ as $\xi$ approaches $\xi_{i+1}$ from left. Denote by $\mathrm{PC}([a,b], \mathbb R^n)$ and $\mathrm{Lip}([a,b], \mathbb R^n)$ the linear spaces of piecewise continuous and Lipschitz continuous functions $x(\cdot) \colon [a,b] \mapsto \mathbb R^n$.

Let $t_0 < \vartheta$ and $h > 0$. Let us denote
$$
\mathrm{PC} = \mathrm{PC}([-h,0],\mathbb R^n),\quad \mathbb G = [t_0,\vartheta] \times \mathbb R^n \times \mathrm{PC}.
$$
Define the following norms on the space $\mathrm{PC}$:
\begin{equation*}
\|w(\cdot)\|_1 = \int\limits_{-h}^0 \|w(\xi)\|\mathrm{d}\xi,\quad \|w(\cdot)\|_\infty = \sup\limits_{\xi\in [-h,0]} \|w(\xi)\|,\quad w(\cdot) \in \mathrm{PC}.
\end{equation*}

Consider a dynamical system which motion is described by the following delay differential equation:
\begin{equation}\label{system}
\dot{x}(\tau) = f(\tau,x(\tau),x(\tau-h),u(\tau)),\quad \tau \in [t_0,\vartheta],\quad x(\tau) \in \mathbb R^n,\quad u(\tau) \in \mathbb U.
\end{equation}
Here, $\tau$ is the time variable, $x(\tau)$ is the state vector at the time $\tau$, $\dot{x}(\tau) = \mathrm{d} x(\tau)/ \mathrm{d} \tau$, $u(\tau)$ is the current control action, $\mathbb U \subset \mathbb R^m$ is a compact set.

Let $(t,z,w(\cdot)) \in \mathbb G$. 
Define
\begin{equation*}
\begin{array}{rl}
\Lambda(t,z,w(\cdot)) = \big\{x(\cdot) \in \mathrm{PC}([t-h,\vartheta],\mathbb R^n)\colon x(\tau) = w(\tau - t),\ \tau \in [t-h,t),\\[0.2cm]
x(\tau) = y(\tau),\ \tau \in [t,\vartheta],\ y(\cdot) \in \mathrm{Lip}([t,\vartheta],\mathbb R^n),\ y(t) = z\big\}.
\end{array}
\end{equation*}
Denote by $\mathfrak{U}(t)$ the set of measurable functions $u(\cdot) \colon [t,\vartheta] \mapsto \mathbb U$. Let $u(\cdot) \in \mathfrak{U}(t)$. By a motion $x(\cdot) = x(\cdot\,|\,t,z,w(\cdot),u(\cdot))$ of system (\ref{system}), we mean a function $x(\cdot) \in \Lambda(t,z,w(\cdot))$ that satisfies equation (\ref{system}) for almost every $\tau \in [t,\vartheta]$.

Consider the following optimal control problem: for each $(t,z,w(\cdot)) \in \mathbb G$, minimize the Bolza cost functional
\begin{equation}\label{J}
J(t,z,w(\cdot),u(\cdot)) = \sigma(x(\vartheta),x_\vartheta(\cdot)) + \int\limits_t^\vartheta f^0(\xi,x(\xi),x(\xi-h),u(\xi))\mathrm{d}\xi,
\end{equation}
over all $u(\cdot) \in \mathfrak{U}(t)$, where $x(\cdot) = x(\cdot\,|\, t,z,w(\cdot),u(\cdot))$ is the motion of system (\ref{system}), $x_\vartheta(\cdot)$ is the function defined by  $x_\vartheta(\xi) = x(\vartheta + \xi)$, $\xi \in [-h,0]$.

We assume that the following conditions hold:
\begin{description}
\item[$(f_1)$] The functions $f(t,x,y,u) \in \mathbb R^n$, $f^0(t,x,y,u) \in \mathbb R$, $t \in [t_0,\vartheta]$, $x,y \in \mathbb R^n$, $u \in \mathbb U$ are continuous.

\item[$(f_2)$] For every $\alpha > 0$, there exists a number $\lambda_f = \lambda_f(\alpha) > 0$ such that
\begin{equation*}
\begin{array}{c}
\big\|f(t,x,y,u) - f(t,x',y',u)\big\| + \big|f^0(t,x,y,u) - f^0(t,x',y',u)\big| \\[0.3cm]
\leq \lambda_f \big(\|x - x'\| + \|y - y'\|\big)
\end{array}
\end{equation*}
for any $t \in [t_0,\vartheta]$, $x,y,x',y' \in O(\alpha) = \{x \in \mathbb R^n \colon \|x\| \leq \alpha\}$ and $u \in \mathbb U$.

\item[$(f_3)$] There exists a constant $c_f > 0$ such that
\begin{equation*}
\big\|f(t,x,y,u)\big\| + \big|f^0(t,x,y,u)\big|\leq c_f (1 + \|x\| + \|y\|)
\end{equation*}
for any $t \in [t_0,\vartheta]$, $x,y \in \mathbb R^n$ and $u \in \mathbb U$.

\item[($\sigma$)] For every $\alpha > 0$, there exists a number $\lambda_\sigma = \lambda_\sigma(\alpha) > 0$ such that
\begin{equation*}
\big|\sigma(z,w(\cdot)) - \sigma(z',w'(\cdot))\big| \leq \lambda_\sigma \big(\|z - z'\| + \|w(\cdot) - w'(\cdot)\|_1\big)
\end{equation*}
for any $(z,w(\cdot)),(z',w'(\cdot)) \in P(\alpha)$, where
\begin{equation*}
P(\alpha) = \big\{(z,w(\cdot)) \in \mathbb R^n \times \mathrm{PC} \colon \|z\| \leq \alpha,\, \|w(\cdot)\|_\infty \leq \alpha\big\}.
\end{equation*}
\end{description}

It is known that, under such conditions, for each $(t,z,w(\cdot)) \in \mathbb G$ and $u(\cdot) \in \mathfrak{U}(t)$, there exists a unique motion $x(\cdot) = x(\cdot\,|\,t,z,w(\cdot),u(\cdot))$ of system (\ref{system}).

The value functional $\rho \colon \mathbb G \mapsto \mathbb R$ in optimal control problem (\ref{system}), (\ref{J}) is defined by
\begin{equation}\label{rho}
\rho(t,z,w(\cdot)) = \inf\limits_{u(\cdot) \in \mathfrak{U}(t)} J(t,z,w(\cdot),u(\cdot)),\quad (t,z,w(\cdot)) \in \mathbb G.
\end{equation}
One can show (following, e.g., the scheme from \cite[p. 553]{Evans_1998}) that, for every $(t,z,w(\cdot)) \in \mathbb G$ and $\tau \in [t,\vartheta]$, the functional $\rho$ satisfies the following equation (a dynamic programming principle):
\begin{equation}\label{dynamic_programming}
\rho(t,z,w(\cdot)) = \inf\limits_{u(\cdot) \in \mathfrak{U}(t)} \bigg(\rho(\tau,x(\tau),x_\tau(\cdot)) + \int\limits_t^\tau f^0(\xi,x(\xi),x(\xi-h),u(\xi))\mathrm{d}\xi\bigg),
\end{equation}
where $x(\cdot) = x(\cdot\,|\, t,z,w(\cdot),u(\cdot))$ is the motion of system (\ref{system}).

In order to obtain a Hamilton-Jacobi-Bellman (HJB) equation as infinitesimal form of equation (\ref{dynamic_programming}), we will use the following definition of differentiability of functionals. Following \cite{Kim_1999,Lukoyanov_2000,Lukoyanov_2001}, a functional $\varphi \colon \mathbb G \mapsto \mathbb R$ is called coinvariantly differentiable (ci-differentiable) at a point $(t,z,w(\cdot)) \in \mathbb G$, $t < \vartheta$ if there exist a number $\partial^{ci}_{t,w}\varphi(t,z,w(\cdot)) \in \mathbb R$ and a vector $\nabla_z\varphi(t,z,w(\cdot)) \in \mathbb R^n$ such that, for any $v \in \mathbb R^n$, $x(\cdot) \in \Lambda(t,z,w(\cdot))$ and $\tau \in [t,\vartheta]$, the following relation holds:
\begin{equation}\label{ci-deriv}
\begin{array}{c}
\varphi(\tau,v,x_\tau(\cdot)) - \varphi(t,z,w(\cdot))  =  \partial^{ci}_{t,w}\varphi(t,z,w(\cdot)) (\tau - t) \\[0.3cm]
+ \langle v - z, \nabla_z \varphi(t,z,w(\cdot)) \rangle + o(|\tau - t| + \|v - z\|),
\end{array}
\end{equation}
where the function $x_\tau(\cdot) \in \mathrm{PC}$ is defined by $x_\tau(\xi)=x(\tau + \xi)$, $\xi \in [-h,0]$, the value $o(\cdot)$ depends on the triplet $\{t,z,x(\cdot)\}$, and $o(\delta)/\delta \to 0$ as $\delta \to +0$. Then $\partial^{ci}_{t,w}\varphi(t,z,w(\cdot))$ is called the ci-derivative of $\varphi$ with respect to $\{t,w(\cdot)\}$ and $\nabla_z \varphi(t,z,w(\cdot))$ is the gradient of $\varphi$ with respect to $z$. Let us note that if $\varphi$ does not depend on the functional variable $w(\cdot)$, then the definition of ci-differentiability coincides with the definition of differentiability of functions.

Define the Hamiltonian of problem (\ref{system}), (\ref{J}) by
\begin{equation}\label{H}
\begin{array}{c}
H(t,x,y,s) = \min\limits_{u \in \mathbb U} \big(\langle f(t,x,y,u),s\rangle + f^0(t,x,y,u)\big),\\[0.3cm]
t \in [t_0,\vartheta],\quad x,y,s \in \mathbb R^n.
\end{array}
\end{equation}
Consider the following Cauchy problem for the HJB equation
\begin{equation}\label{HJ}
\begin{array}{c}
\partial^{ci}_{t,w} \varphi(t,z,w(\cdot)) + H(t,z,w(-h),\nabla_z \varphi(t,z,w(\cdot))) = 0,\\[0.3cm]
(t,z,w(\cdot)) \in \mathbb G,\quad t < \vartheta,
\end{array}
\end{equation}
and the terminal condition
\begin{equation}\label{HJ0}
\varphi(\vartheta,z,w(\cdot)) =\sigma(z,w(\cdot)),\quad (\vartheta,z,w(\cdot)) \in \mathbb G.
\end{equation}
Define the class of functionals in which we will search a solution of this problem. Denote by $\Phi$ the set of functionals $\varphi = \varphi(t,z,w(\cdot)) \in \mathbb R$, $(t,z,w(\cdot)) \in \mathbb G$ which are continuous with respect to $t$ and satisfy the following Lipschitz condition: for every $\alpha > 0$, there exists a number $\lambda_\varphi = \lambda_\varphi(\alpha) > 0$ such that
\begin{equation}\label{Phi_lip}
|\varphi(t,z,w(\cdot)) - \varphi(t,z',w'(\cdot))| \leq \lambda_\varphi\big(\|z - z'\| + \|w(\cdot) - w'(\cdot)\|_1\big)
\end{equation}
for any $t \in [t_0,\vartheta]$ and $(z,w(\cdot)), (z',w'(\cdot)) \in P(\alpha)$. The choice of this class is motivated, in particular, the inclusion $\rho \in \Phi$, which will be shown in Lemma~\ref{lem_rho_Phi}.

The following theorem establishes the relation between problem (\ref{HJ}), (\ref{HJ0}) and the value functional $\rho$ in the case when $\rho$ is ci-differentiable.

\begin{theorem}\label{thm1}
The following statements hold:
\begin{itemize}
\item[1.] If a functional $\varphi \in \Phi$ is ci-differentiable at each point $(t,z,w(\cdot)) \in \mathbb G$, $t < \vartheta$, satisfies HJB equation {\rm (\ref{HJ})} at these points and satisfies terminal condition {\rm (\ref{HJ0})}, then the identity $\varphi \equiv \rho$ holds.

\item[2.] If the value functional $\rho$ is ci-differentiable at a point $(t,z,w(\cdot)) \in \mathbb G$, $t < \vartheta$, then it satisfies HJB equation {\rm (\ref{HJ})} at this point.
\end{itemize}
\end{theorem}
The proof of this theorem is described after the Theorem \ref{thm2}.

For the case when $\rho$ is not ci-differentiable, definitions of generalized (minimax and viscosity) solutions of problem (\ref{HJ}), (\ref{HJ0}) are given below.

Taking the constant $c_f > 0$ form $(f_3)$, we denote
\begin{equation}\label{F}
F(x,y) = \big\{l \in \mathbb R^n\colon \|l\| \leq c_f (1 + \|x\| + \|y\|) \big\} \subset \mathbb R^n,\quad x,y \in \mathbb R^n.
\end{equation}
Let $(t,z,w(\cdot)) \in \mathbb G$. Denote by $X(t,z,w(\cdot))$ the set of the functions $x(\cdot) \in \Lambda(t,z,w(\cdot))$ that satisfy the following delay differential inclusion:
\begin{equation}\label{inclution}
\dot{x}(\tau) \in F(x(\tau),x(\tau-h))\text{ for a.e. } \tau \in [t,\vartheta].
\end{equation}
Note that the set $X(t,z,w(\cdot))$ is not empty. In particular, for each $u(\cdot) \in \mathfrak{U}(t)$, the motion $x(\cdot) = x(\cdot\,|\,t,z,w(\cdot),u(\cdot))$ of system (\ref{system}) satisfies the inclusion
\begin{equation}\label{xX}
x(\cdot) \in X(t,z,w(\cdot)).
\end{equation}

\begin{definition}
A functional $\varphi \colon \mathbb G \mapsto \mathbb R$ is called a minimax solution of problem {\rm (\ref{HJ}), (\ref{HJ0})} if $\varphi$ satisfies the inclusion $\varphi \in \Phi$, terminal condition {\rm (\ref{HJ0})} and the following inequalities:
\begin{subequations}\label{minmax_sol}
\begin{align}
\begin{array}{rl}
\inf\limits_{x(\cdot) \in X(t,z,w(\cdot))} & \bigg(\varphi(\tau,x(\tau),x_\tau(\cdot)) - \varphi(t,z,w(\cdot)) \\[0.2cm]
& \displaystyle\quad + \int\limits_t^\tau\big(H(\xi,x(\xi),x(\xi-h),s) - \langle \dot{x}(\xi),s \rangle\big)\mathrm{d}\xi\bigg) \leq 0, \label{minmax_sol:v}
\end{array} \\[0.1cm]
\begin{array}{rl}
\sup\limits_{x(\cdot) \in X(t,z,w(\cdot))} & \bigg(\varphi(\tau,x(\tau),x_\tau(\cdot)) - \varphi(t,z,w(\cdot)) \\[0.2cm]
& \displaystyle\quad + \int\limits_t^\tau\big(H(\xi,x(\xi),x(\xi-h),s) - \langle \dot{x}(\xi),s \rangle\big)\mathrm{d}\xi\bigg) \geq 0, \label{minmax_sol:n}
\end{array}
\end{align}
\end{subequations}
for any $(t,z,w(\cdot)) \in \mathbb G$, $t < \vartheta$, $\tau \in (t,\vartheta]$ and $s \in \mathbb R^n$.
\end{definition}

By analogy with \cite{Lukoyanov_2010a}, lower and upper right directional derivatives of a functional $\varphi \colon \mathbb G \mapsto \mathbb R$ along $l \in \mathbb R^n$ at $(t,z,w(\cdot)) \in \mathbb G$, $t < \vartheta$ are defined by
\begin{subequations}\label{direct_deriv}
\begin{align}
\partial^-_{l} \varphi(t,z,w(\cdot)) &
= \liminf\limits_{\tau \to t + 0} \displaystyle\frac{\varphi(\tau,x^l(\tau),x^l_\tau(\cdot)) - \varphi(t,z,w(\cdot))}{\tau - t},\label{direct_deriv:-} \\[0.0cm]
\partial^+_{l} \varphi(t,z,w(\cdot)) &
= \limsup\limits_{\tau \to t + 0} \displaystyle\frac{\varphi(\tau,x^l(\tau),x^l_\tau(\cdot)) - \varphi(t,z,w(\cdot))}{\tau - t},\label{direct_deriv:+}
\end{align}
\end{subequations}
where $x^l(\cdot) \in \Lambda(t,z,w(\cdot))$ and $x^l(\tau) = z + l (\tau - t)$, $\tau \in [t,\vartheta]$.

The following sets are called the subdifferential $D^-\varphi(t,z,w(\cdot))$ and the superdifferential $D^+\varphi(t,z,w(\cdot))$ of the functional $\varphi$ at $(t,z,w(\cdot)) \in \mathbb G$, $t < \vartheta$:
\begin{subequations}\label{D}
\begin{align}
\begin{array}{rl}
D^-\varphi(t,z,w(\cdot)) = & \big\{(p_0,p) \in \mathbb R \times \mathbb R^n\colon  \\[0.1cm]
& \quad p_0 + \langle l,p \rangle \leq \partial^-_l \varphi(t,z,w(\cdot)),\, l \in \mathbb R^n\big\},\label{D:-}
\end{array} \\[0.1cm]
\begin{array}{rl}
D^+\varphi(t,z,w(\cdot)) = & \big\{(q_0,q) \in \mathbb R \times \mathbb R^n\colon  \\[0.1cm]
& \quad q_0 + \langle l,q \rangle \geq \partial^+_l \varphi(t,z,w(\cdot)),\, l \in \mathbb R^n\big\}.\label{D:+}
\end{array}
\end{align}
\end{subequations}
Note that if a functional $\varphi$ is ci-differentiable at $(t,z,w(\cdot)) \in \mathbb G$, $t < \vartheta$, then
\begin{eqnarray}\label{ci-phi}
& \hspace{-0.3cm} \partial^-_l \varphi(t,z,w(\cdot)) = \partial^+_l \varphi(t,z,w(\cdot)) = \partial^{ci}_{t,w}\varphi(t,z,w(\cdot)) + \langle l,\nabla_z\varphi(t,z,w(\cdot)) \rangle,&\nonumber \\[0.2cm]
& \hspace{-0.3cm} D^-\varphi(t,z,w(\cdot)) = \big\{(p_0,p)\colon p_0 \leq \partial^{ci}_{t,w}\varphi(t,z,w(\cdot)),\, p = \nabla_z\varphi(t,z,w(\cdot))\big\},&\\[0.2cm]
& \hspace{-0.3cm} D^+\varphi(t,z,w(\cdot)) = \big\{(q_0,q)\colon q_0 \geq \partial^{ci}_{t,w}\varphi(t,z,w(\cdot)),\, q = \nabla_z\varphi(t,z,w(\cdot))\big\}.&\nonumber
\end{eqnarray}

\begin{definition}
A functional $\varphi \colon \mathbb G \mapsto \mathbb R$ is called a viscosity solution of problem {\rm (\ref{HJ}), (\ref{HJ0})} if $\varphi$ satisfies the inclusion $\varphi \in \Phi$, terminal condition {\rm (\ref{HJ0})} and the following inequalities:
\begin{subequations}\label{visc_sol}
\begin{align}
p_0 + H(t,z,w(-h),p) & \leq 0,\quad (p_0,p) \in D^-\varphi(t,z,w(\cdot)),\label{visc_sol:-} \\[0.2cm]
q_0 + H(t,z,w(-h),q) & \geq 0,\quad (q_0,q) \in D^+\varphi(t,z,w(\cdot))\label{visc_sol:+}
\end{align}
\end{subequations}
for any $(t,z,w(\cdot)) \in \mathbb G$, $t < \vartheta$.
\end{definition}

\begin{theorem}\label{thm2}
For $\varphi \colon \mathbb G \mapsto \mathbb R$, the following statements are equivalent:
\begin{description}
\item[(a)] The identity $\varphi \equiv \rho$ holds.
\vspace{0.1cm}

\item[(b)] $\varphi$ is a minimax solution of problem {\rm (\ref{HJ}), (\ref{HJ0})}.
\vspace{0.1cm}

\item[(c)] $\varphi \in \Phi$ satisfies terminal condition {\rm (\ref{HJ0})} and the following inequalities:
\begin{subequations}\label{ineq_for_deriv}
\begin{align}
\inf\limits_{l \in F(z,w(-h))} \big(\partial^-_l \varphi(t,z,w(\cdot)) + H(t,z,w(-h),s) - \langle l,s \rangle\big) & \leq 0,\label{ineq_for_deriv:-} \\[0.0cm]
\sup\limits_{l \in F(z,w(-h))} \big(\partial^+_l \varphi(t,z,w(\cdot)) + H(t,z,w(-h),s) - \langle l,s \rangle\big) & \geq 0\label{ineq_for_deriv:+}
\end{align}
\end{subequations}
for any $(t,z,w(\cdot)) \in \mathbb G$, $t < \vartheta$ and $s \in \mathbb R^n$.
\vspace{0.1cm}

\item[(d)] $\varphi$ is a viscosity solution of problem {\rm (\ref{HJ}), (\ref{HJ0})}.
\end{description}
\end{theorem}
In particular, this theorem establishes the existence and uniqueness of the minimax and viscosity solutions since the value functional $\rho$ is uniquely defined.

Note that Theorem \ref{thm1} follows from the equivalence of statements $(a)$ and $(d)$ if we take into account (\ref{ci-phi}). Below in the paper, auxiliary properties of system (\ref{system}) and inclusion (\ref{inclution}) will be given and Theorem \ref{thm2} will be proved.

\section{Properties of Time-Delay Systems}

\begin{proposition}\label{ass_X}
For every $\alpha > 0$, there exist numbers $\alpha_X = \alpha_X(\alpha) > \alpha$ and $\lambda_X = \lambda_X(\alpha) > 0$ such that
\begin{equation}\label{ass_X-0}
(x(\tau),x_\tau(\cdot)) \in P(\alpha_X),\quad \|x(\tau) - x(\tau')\| \leq \lambda_X |\tau - \tau'|,\quad \tau,\tau' \in [t,\vartheta]
\end{equation}
for each $t \in [t_0,\vartheta]$, $(z,w(\cdot)) \in P(\alpha)$ and $x(\cdot) \in X(t,z,w(\cdot))$.
\end{proposition}
\proof Let $\alpha > 0$. Put $\alpha_* = (1 + c_f h) \alpha + c_f(\vartheta - t_0)$, $\alpha_X = \alpha_* e^{2c_f(\vartheta - t_0)}$ and $\lambda_X = c_f(1 + 2\alpha_X)$.
Let $t \in [t_0,\vartheta]$, $(z,w(\cdot)) \in P(\alpha)$ and $x(\cdot) \in X(t,z,w(\cdot))$. Then, according to (\ref{F}), (\ref{inclution}), we derive
\begin{equation*}
\|x(\tau)\| \leq \|z\| + c_f \int\limits_t^\tau \big(1 + \|x(\xi)\| + \|x(\xi - h)\|\big) \mathrm{d}\xi \leq \alpha_* + 2 c_f \int\limits_t^\tau \|x(\xi)\|\mathrm{d}\xi.
\end{equation*}
Therefore, applying Bellman-Gronwall lemma (see, e.g., \cite[p. 31]{Bellman_Cooke_1963}), we obtain $(x(\tau),x_\tau(\cdot)) \in P(\alpha_X)$, $\tau \in [t,\vartheta]$. Then, from (\ref{inclution}), we deduce $\|\dot{x}(\tau)\| \leq \lambda_X$ for almost every $\tau \in [t,\vartheta]$, which concludes the proof.\hfill $\square$

For $(t,z,w(\cdot)) \in \mathbb G$ and $u(\cdot) \in \mathfrak{U}(t)$, we denote
\begin{equation}\label{I}
I(t,z,w(\cdot),u(\cdot)) = \int\limits_t^\vartheta f^0(\xi,x(\xi),x(\xi-h),u(\xi))\mathrm{d}\xi,
\end{equation}
where $x(\cdot) = x(\cdot\,|\,t,z,w(\cdot),u(\cdot))$ is the motion of system (\ref{system}).

\begin{proposition}\label{ass_Phi_lip}
For every $\alpha > 0$, there exists a number $\lambda_* = \lambda_*(\alpha) > 0$ such that, for each $t \in [t_0,\vartheta]$, $(z,w(\cdot)), (z',w'(\cdot)) \in P(\alpha)$ and $u(\cdot) \in \mathfrak{U}(t)$, the motions $x(\cdot) = x(\cdot\,|\,t,z,w(\cdot),u(\cdot))$ and $x'(\cdot) = x(\cdot\,|\,t,z',w'(\cdot),u(\cdot))$ of system {\rm (\ref{system})} satisfy the inequality
\begin{equation*}
\|x(\vartheta) - x'(\vartheta)\| + \|x_\vartheta(\cdot) - x'_\vartheta(\cdot)\|_1 + |I(t,z,w(\cdot),u(\cdot)) - I(t,z',w'(\cdot),u(\cdot))|
\end{equation*}
\begin{equation}\label{ass_Phi_lip_0}
\leq \lambda_* \big(\|z - z'\| + \|w(\cdot) - w'(\cdot)\|_1\big).
\end{equation}
\end{proposition}
\proof Let $\alpha > 0$. According to Proposition \ref{ass_X} and condition $(f_2)$, let us define the numbers $\alpha_X = \alpha_X(\alpha) > \alpha$ and $\lambda_f = \lambda_f(\alpha_X) > 0$.
Put
\begin{equation}\label{ass_Phi_lip-lambda_s}
\lambda_* = (1 + \lambda_f)\big(2 + (1 + 2\lambda_f)(\vartheta - t_0) e^{2\lambda_f(\vartheta - t_0)}\big).
\end{equation}
Let $t \in [t_0,\vartheta]$, $(z,w(\cdot)), (z',w'(\cdot)) \in P(\alpha)$ and $u(\cdot) \in \mathfrak{U}(t)$. Then, in accordance with (\ref{xX}), for the motions $x(\cdot) = x(\cdot\,|\,t,z,w(\cdot),u(\cdot))$ and $x'(\cdot) = x(\cdot\,|\,t,z',w'(\cdot),u(\cdot))$, we have
\begin{equation}\label{ass_Phi_lip-2}
\begin{array}{c}
\displaystyle \|x(\tau) - x'(\tau)\| + |I(t,z,w(\cdot),u(\cdot)) - I(t,z',w'(\cdot),u(\cdot))|\\[0.1cm]
\displaystyle \leq \|z - z'\| + \lambda_f \|w(\cdot) - w'(\cdot)\|_1 + 2 \lambda_f \int\limits_t^\tau \|x(\xi) - x'(\xi)\|\mathrm{d}\xi,\quad \tau \in [t,\vartheta].
\end{array}
\end{equation}
Then, applying Bellman-Gronwall lemma (see, e.g., \cite[p. 31]{Bellman_Cooke_1963}), we obtain
\begin{equation}\label{ass_Phi_lip-3}
\|x(\tau) - x'(\tau)\| \leq \big(\|z - z'\| + \lambda_f \|w(\cdot) - w'(\cdot)\|_1\big) e^{2 \lambda_f (\vartheta - t_0)},\quad \tau \in [t,\vartheta].
\end{equation}
From (\ref{ass_Phi_lip-lambda_s})--(\ref{ass_Phi_lip-3}), taking into account the estimate
\begin{equation*}
\|x_\tau(\cdot) - x'_\tau(\cdot)\|_1 \leq \|w(\cdot) - w'(\cdot)\|_1 + (\vartheta - t_0) \max\limits_{\xi \in [t,\tau]}\|x(\xi) - x'(\xi)\|,\quad
\tau \in [t,\vartheta],
\end{equation*}
we conclude (\ref{ass_Phi_lip_0}). \hfill $\square$

\begin{proposition}\label{ass_w}
Let $(t,z,w(\cdot)) \in \mathbb G$. For every $\varepsilon > 0$, there exists a number $\delta = \delta(\varepsilon) > 0$ such that, for every $t' \in [t,t + \delta]$ and $x(\cdot) \in X(t,z,w(\cdot))$, the inequality $\|w(\cdot) - x_{t'}(\cdot)\|_1 \leq \varepsilon$ holds.
\end{proposition}
This proposition one can proved, using approximation of $w(\cdot)$ by a Lipschitz function (see, e.g., \cite[p. 214]{Natanson_1960}) and Proposition \ref{ass_X}.

\begin{proposition}\label{ass_Phi_coun}
For every $(z,w(\cdot)) \in \mathbb R^n \times \mathrm{PC}$ and every $\varepsilon > 0$, there exists a number $\delta_* = \delta_*(z,w(\cdot),\varepsilon) > 0$ such that, for every $t,t' \in [t_0,\vartheta]$ with $|t - t'| \leq \delta_*$ and every $u(\cdot) \in \mathfrak{U}(t)$, there exists a function $u'(\cdot) \in \mathfrak{U}(t')$ such that the motions $x(\cdot) = x(\cdot\,|\,t,z,w(\cdot),u(\cdot))$ and $x'(\cdot) = x(\cdot\,|\,t',z,w(\cdot),u'(\cdot))$ of system {\rm (\ref{system})} satisfy the inequality
\begin{equation*}
\|x(\vartheta) - x'(\vartheta)\| + \|x_\vartheta(\cdot) - x'_\vartheta(\cdot)\|_1 + |I(t,z,w(\cdot),u(\cdot)) - I(t',z,w(\cdot),u'(\cdot))| \leq \varepsilon.
\end{equation*}
\end{proposition}
\proof Let $(z,w(\cdot)) \in \mathbb R^n \times \mathrm{PC}$ and $\varepsilon > 0$. Take a number $\alpha > 0$ so that $(z,w(\cdot)) \in P(\alpha)$ (for instance, $\alpha = \max\{\|z\|,\|w(\cdot)\|_\infty\} + 1$). In accordance with Proposition \ref{ass_X}, let us the take numbers $\alpha_X = \alpha_X(\alpha) > \alpha$ and $\lambda_X = \lambda_X(\alpha) > 0$. Then, taking into account $(f_3)$, for any $t \in [t_0,\vartheta]$, $t' \in [t,\vartheta]$ and $u(\cdot) \in \mathfrak{U}(t)$, the motion $x(\cdot) = x(\cdot\,|\,t,z,w(\cdot),u(\cdot))$ satisfies
\begin{equation}\label{ass_Phi_coun-1}
\begin{array}{c}
\|x(t') - z\| \leq \lambda_X (t' - t), \\[0.2cm]
|I(t',x(t'),x_{t'}(\cdot),u(\cdot)) - I(t,z,w(\cdot),u(\cdot))| \leq c_f (1 + 2 \alpha_X) (t' - t).
\end{array}
\end{equation}
Doe to Proposition \ref{ass_Phi_lip}, let us take $\lambda_* = \lambda_*(\alpha_X) > 0$. Let $t,t' \in [t_0,\vartheta]$ and $u(\cdot) \in \mathfrak{U}(t)$. If $t' \geq t$ then, defining $u'(\cdot) \in \mathfrak{U}(t')$ by $u'(\tau) = u(\tau)$, $\tau \in [t',\vartheta]$, for $x(\cdot) = x(\cdot\,|\,t,z,w(\cdot),u(\cdot))$ and $x'(\cdot) = x(\cdot\,|\,t',z,w(\cdot),u'(\cdot))$, we have
\begin{equation*}
\|x(\vartheta) - x'(\vartheta)\| + \|x_\vartheta(\cdot) - x'_\vartheta(\cdot)\|_1 + |I(t',x(t'),x_{t'}(\cdot),u(\cdot)) - I(t',z,w(\cdot),u'(\cdot))|
\end{equation*}
\begin{equation}\label{ass_Phi_coun-2}
\leq \lambda_* \big(\|x(t') - z\| + \|x_{t'}(\cdot) - w(\cdot)\|_1\big),
\end{equation}
If $t' < t$, then, defining $u'(\cdot) \in \mathfrak{U}(t')$ by $u'(\tau) = u(t)$, $\tau \in [t',t)$ and $u'(\tau) = u(\tau)$, $\tau \in [t,\vartheta]$, for $x(\cdot) = x(\cdot\,|\,t,z,w(\cdot),u(\cdot))$ and $x'(\cdot) = x(\cdot\,|\,t',z,w(\cdot),u'(\cdot))$, we have
\begin{equation*}
\|x(\vartheta) - x'(\vartheta)\| + \|x_\vartheta(\cdot) - x'_\vartheta(\cdot)\|_1 + |I(t,z,w(\cdot),u(\cdot)) - I(t,x'(t),x'_t(\cdot),u'(\cdot))|
\end{equation*}
\begin{equation}\label{ass_Phi_coun-3}
\leq \lambda_* \big(\|z - x'(t)\| + \|w(\cdot) - x'_t(\cdot)\|_1\big).
\end{equation}
The proof of this proposition follows from (\ref{ass_Phi_coun-1})--(\ref{ass_Phi_coun-3}) and Proposition \ref{ass_w}.\hfill $\square$

\section{Proof of Theorem \ref{thm2}}

\subsection{Proof $(a) \Rightarrow (d)$}

Notes that, under conditions $(f_1)$--$(f_3)$, the following conditions hold:
\begin{description}
\item[$(H_1)$] The function $H(t,x,y,s) \in \mathbb R^n$, $t \in [t_0,\vartheta]$, $x,y,s \in \mathbb R^n$ defined by (\ref{H}) is continuous.

\item[$(H_2)$] For every $\alpha > 0$, there exists a number $\lambda_H = \lambda_H(\alpha) > 0$ such that
\begin{equation*}
|H(t,x,y,s) - H(t,x',y',s)| \leq \lambda_H (\|x - x'\| + \|y - y'\|) (1 + \|s\|)
\end{equation*}
for any $t \in [t_0,\vartheta]$, $x,y,x',y' \in O(\alpha)$ and $u \in \mathbb U$.

\item[$(H_3)$] For any $t \in [t,\vartheta]$ and $x,y,s \in \mathbb R^n$ the following equalities hold:
\begin{equation*}
\begin{array}{rcl}
H(t,x,y,s) & = & \max\limits_{q \in \mathbb R^n} \min\limits_{f \in F(x,y)} \big(H(t,x,y,q) + \langle f,s-q \rangle\big)\\[0.3cm]
& = & \min\limits_{p \in \mathbb R^n} \max\limits_{f \in F(x,y)} \big(H(t,x,y,p) + \langle f,s-p \rangle\big).
\end{array}
\end{equation*}

\item[($F$)] The multi-valued mapping $F(x,y) \subset \mathbb R^n$, $x,y \in \mathbb R^n$ defined by (\ref{F}) is continuous with respect to Hausdorff metric.
\end{description}

Conditions $(H_1)$--$(H_3)$ can be proved by analogy with \cite[p. 14]{Subbotin_1995} and \cite[p. 177]{Subbotin_1995}. Condition $(F)$ directly follows from (\ref{F}).

\begin{lemma}\label{lem_rho_Phi}
The value functional $\rho$ satisfies the inclusion $\rho \in \Phi$.
\end{lemma}
\proof Let us first prove that $\rho$ is uniformly continuous with respect to $t$. Let $(z,w(\cdot)) \in \mathbb R^n \times \mathrm{PC}$ and $\varepsilon > 0$. Take a number $\alpha > 0$ so that $(z,w(\cdot)) \in P(\alpha)$. In accordance with Proposition~\ref{ass_X} and condition $(\sigma)$, define the numbers $\alpha_X = \alpha_X(\alpha) > \alpha$ and $\lambda_\sigma = \lambda_\sigma(\alpha_X)$. Then, according to Proposition~\ref{ass_Phi_coun}, there exists $\delta_* = \delta_*(\varepsilon/(2\lambda_\sigma + 2))$ such that, for every $t,t' \in [t_0,\vartheta]$ with $|t - t'| \leq \delta_*$ and every $u(\cdot) \in \mathfrak{U}(t)$, there exists $u'(\cdot) \in \mathfrak{U}(t')$ such that
\begin{equation}\label{ass_rho_Phi-1}
|J(t,z,w(\cdot),u(\cdot)) - J(t',z,w(\cdot),u'(\cdot))| \leq \varepsilon/2,
\end{equation}
where $J(t,z,w(\cdot),u(\cdot))$, $J(t',z,w(\cdot),u'(\cdot))$ are defined by (\ref{J}), (\ref{I}).

Let $t,t' \in [t_0,\vartheta]$, $|t - t'| \leq \delta_*$. Due to (\ref{rho}), there exists $u(\cdot) \in \mathfrak{U}(t)$ such that
\begin{equation}\label{ass_rho_Phi-2}
J(t,z,w(\cdot),u(\cdot)) \leq \rho(t,z,w(\cdot)) + \varepsilon/2.
\end{equation}
Defining $u'(\cdot) \in \mathfrak{U}(t')$ so that (\ref{ass_rho_Phi-1}) holds and taking into account (\ref{rho}), we obtain
\begin{equation*}
\rho(t',z,w(\cdot)) - \rho(t,z,w(\cdot)) \leq J(t',z,w(\cdot),u'(\cdot)) - J(t,z,w(\cdot),u(\cdot)) + \varepsilon/2 \leq \varepsilon.
\end{equation*}
Thus, uniform continuity of $\rho$ with respect to $t$ has been proved. Also, using the similar way of proving, basing on Proposition \ref{ass_Phi_lip}, one can show that $\rho$ satisfies Lipschitz condition (\ref{Phi_lip}).
%
\hfill $\square$

\begin{lemma}\label{lem_partial_k}
Let $(t,z,w(\cdot)) \in \mathbb G$, $t < \vartheta$ and $\tau_k \in [t,\vartheta]$, $x^{(k)}(\cdot) \in X(t,z,w(\cdot))$, $k \in \mathbb N$. Let $\tau_k \to t$ and $l_k = (x^{(k)}(\tau_k) - z) / (\tau_k - t) \to l_* \in \mathbb R^n$ as $k \to \infty$. Then the following inequalities hold:
\begin{subequations}\label{lem_partial_k-0}
\begin{align}
\displaystyle\partial^-_{l_*}\rho(t,z,w(\cdot))
\leq \liminf\limits_{k \to \infty} \frac{\rho(\tau_k,x^{(k)}(\tau_k),x^{(k)}_{\tau_k}(\cdot)) - \rho(t,z,w(\cdot))}{\tau_k - t}, \label{lem_partial_k-0:-}\\[0.1cm]
\displaystyle\partial^+_{l_*}\rho(t,z,w(\cdot))
\geq \limsup\limits_{k \to \infty} \frac{\rho(\tau_k,x^{(k)}(\tau_k),x^{(k)}_{\tau_k}(\cdot)) - \rho(t,z,w(\cdot))}{\tau_k - t}. \label{lem_partial_k-0:+}
\end{align}
\end{subequations}
\end{lemma}
\proof Take a number $\alpha > 0$ so that $(z,w(\cdot)) \in P(\alpha)$. In accordance with Proposition \ref{ass_X}, define the number $\lambda_X = \lambda_X(\alpha)$. Then we have $\|l_k\| \leq \lambda_X$, $k \in \mathbb N$ and $\|l_*\| \leq \lambda_X$. Put $\alpha_* = \alpha + \lambda_X(\vartheta - t)$. Then, according to Lemma \ref{lem_rho_Phi}, there exists a number $\lambda_\varphi = \lambda_\varphi(\alpha_*) > 0$ such that
\begin{equation*}
|\rho(\tau_k,x^{(k)}(\tau_k),x^{(k)}_{\tau_k}(\cdot)) - \rho(\tau_k,x^{l_*}(\tau_k),x^{l_*}_{\tau_k}(\cdot))|
\leq \lambda_\varphi \big(\|l_k - l_*\| + \lambda_X (\tau_k - t)\big) (\tau_k - t).
\end{equation*}
Consequently, taking into account (\ref{direct_deriv:-}), (\ref{direct_deriv:+}), we obtain (\ref{lem_partial_k-0:-}), (\ref{lem_partial_k-0:+}).\hfill $\square$

\proof $(a) \Rightarrow (d)$. According to (\ref{rho}), the value functional $\rho$ satisfies terminal condition (\ref{HJ0}). Due to Lemma~\ref{lem_rho_Phi}, the inclusion $\rho \in \Phi$ holds. Thus, for the proof of statement $(d)$, it is sufficient to show that $\rho$ satisfies (\ref{visc_sol:-}), (\ref{visc_sol:+}).

Let $(t,z,w(\cdot)) \in \mathbb G$, $t < \vartheta$. Let $(p_0,p) \in D^-\varphi(t,z,w(\cdot))$ and $(q_0,q) \in D^+\varphi(t,z,w(\cdot))$. Since $w(\cdot) \in \mathrm{PC}$, then there exists $\tau_* \in (t,\min\{\vartheta,t+h\})$ such that $w(\cdot)$ is continuous on $[-h,\tau_*-t-h]$. Then, according to conditions $(f_1)$, $(H_1)$, taking into account Proposition~\ref{ass_X}, there exists a sequence $\tau_k \in [t,\tau_*]$, $k \in \mathbb N$ such that $\tau_k \to t$ as $k \to \infty$ and
\begin{equation}\label{1to4-coun}
\begin{array}{rcl}
\big|H(\xi,x(\xi),w(\xi-t-h),p) - H(t,z,w(-h),p)\big| & \leq & 1/k,\\[0.3cm]
\big|f^0(\xi,x(\xi),w(\xi-t-h),u) - f^0(t,z,w(-h),u)\big| & \leq & 1/k,\\[0.3cm]
\big\|f(\xi,x(\xi),w(\xi-t-h),u) - f(t,z,w(-h),u)\big\| & \leq & 1/k
\end{array}
\end{equation}
for any $\xi \in [t,\tau_k]$, $x(\cdot) \in X(t,z,w(\cdot))$, $u \in \mathbb U$.

Let us prove that $\rho$ satisfies (\ref{visc_sol:-}). Due to (\ref{dynamic_programming}), for each $k \in \mathbb N$, there exists a function $u^{(k)}(\cdot) \in \mathfrak{U}(t)$ such that the motion $x^{(k)}(\cdot) = x(\cdot\,|\, t,z,w(\cdot),u^{(k)}(\cdot))$ satisfies the inequality
\begin{equation*}
\begin{array}{c}
\displaystyle\rho(\tau_k,x^{(k)}(\tau_k),x^{(k)}_{\tau_k}(\cdot)) + \int\limits_{t}^{\tau_k} f^0(\xi,x^{(k)}(\xi),w(\xi-t-h),u^{(k)}(\xi))\mathrm{d}\xi \\[0.3cm]
\leq \rho(t,z,w(\cdot)) + (\tau_k - t) / k.
\end{array}
\end{equation*}
From this inequality, using (\ref{H}), (\ref{xX}), (\ref{1to4-coun}), we get
\begin{equation}\label{1to4-1}
\begin{array}{c}
\rho(\tau_k,x^{(k)}(\tau_k),x^{(k)}_{\tau_k}(\cdot)) + H(t,z,w(-h),p) (\tau_k - t) - \langle x^{(k)}(\tau_k) - z,p \rangle \\[0.3cm]
\leq \rho(t,z,w(\cdot)) + 2 (\tau_k - t) / k.
\end{array}
\end{equation}
We denote $l_k = (x^{(k)}(\tau_k) - z)/(\tau_k - t)$. Then, taking into account (\ref{xX}), for the number $\lambda_X > 0$ from Proposition~\ref{ass_X}, we have $\|l_k\| \leq \lambda_X$, $k \in \mathbb N$. Without loss of generality, we can suppose that $l_k \to l_* \in \mathbb R^n$ as $k \to \infty$. Then, applying Lemma \ref{lem_partial_k} and inequality (\ref{1to4-1}), we obtain
\begin{equation*}
\partial^-_{l_*} \rho(t,z,w(\cdot)) \leq - H(t,z,w(-h),p) + \langle l_*,p \rangle.
\end{equation*}
Consequently, in accordance with (\ref{D:-}), we conclude (\ref{visc_sol:-}).

Let us prove that $\rho$ satisfies (\ref{visc_sol:+}). According to (\ref{H}), there exists $u_* \in \mathbb U$ such that
\begin{equation}\label{1to4-u}
\langle f(t,z,w(-h),u_*),q \rangle + f^0(t,z,w(-h),u_*) = H(t,z,w(-h),q).
\end{equation}
Consider the motion $x(\cdot) = x(\cdot\,|\,t,z,w(\cdot),u(\cdot) \equiv u_*)$. Due to (\ref{dynamic_programming}), we have
\begin{equation*}
\rho(\tau,x(\tau),x_\tau(\cdot)) + \int_t^\tau f^0(\xi,x(\xi),w(\xi-t-h),u_*)\mathrm{d}\xi \geq \rho(t,z,w(\cdot)),\quad \tau \in [t,\tau_*].
\end{equation*}
Then, using (\ref{1to4-coun}), (\ref{1to4-u}), we get
\begin{equation*}
\rho(\tau_k,x(\tau_k),x_{\tau_k}(\cdot)) + (\tau_k - t) \big(H(t,z,w(-h),q) - \langle f(t,z,w(-h),u_*),q \rangle\big)
\end{equation*}
\begin{equation}\label{1to4-2}
\geq \rho(t,z,w(\cdot)) - (\tau_k - t) / k.
\end{equation}
We denote $l_* = f(t,z,w(-h),u_*)$ and $l_k = (x(\tau_k) - z)/(\tau_k - t)$, $k \in \mathbb N$. Then, according to (\ref{1to4-coun}), we have $l_k \to l_*$ as $k \to \infty$. Applying Lemma \ref{lem_partial_k} and inequality (\ref{1to4-2}), we obtain
\begin{equation*}
\partial^+_{l_*} \rho(t,z,w(\cdot)) \geq - H(t,z,w(-h),q) + \langle l_*,q \rangle.
\end{equation*}
Consequently, in accordance with (\ref{D:+}), we conclude (\ref{visc_sol:+}).\hfill $\square$

\subsection{Proof $(d) \Rightarrow (c)$}

\begin{lemma}\label{lem_lim}
Let $\varphi \in \Phi$, $(t,z,w(\cdot)) \in \mathbb G$, $\tau_k,\tau_* \in [t,\vartheta]$, $z_k,z_* \in \mathbb R^n$ and $x^{(k)}(\cdot),x^*(\cdot) \in \Lambda(t,z,w(\cdot))$, $k \in \mathbb N$. Let
\begin{equation}\label{lem_upsilon-0}
\tau_k \to \tau_*,\quad
z_k \to z_*,\quad
\max\limits_{\xi \in [t,\vartheta]}\|x^{(k)}(\xi) - x^*(\xi)\| \to 0\quad  \text{ as }\quad k \to \infty.
\end{equation}
Then $\varphi(\tau_k,z_k,x^{(k)}_{\tau_k}(\cdot)) \to \varphi(\tau_*,z_*,x^*_{\tau_*}(\cdot))$ as $k \to \infty$.
\end{lemma}
Using approximation of the function $x^*(\cdot)$ by a Lipschitz function (see, e.g., \cite[p. 214]{Natanson_1960}), we can get that $\|x^*_{\tau_k}(\cdot) - x^*_{\tau_*}(\cdot)\|_1 \to 0$ as $k \to \infty$. Then the lemma follows from (\ref{lem_upsilon-0}), according to the inclusion $\varphi \in \Phi$

\begin{lemma}\label{lem_partial_coun}
Let $\varphi \in \Phi$ and $(t,z,w(\cdot)) \in \mathbb G$, $t < \vartheta$. If there exists a vector $l_0 \in \mathbb R^n$ such that $\partial^-_{l_0}\varphi(t,z,w(\cdot)) \in \mathbb R$, then $\partial^-_l\varphi(t,z,w(\cdot)) \in \mathbb R$ for every $l \in \mathbb R^n$, and the function $\mathbb R^n \ni l \mapsto \partial^-_l\varphi(t,z,w(\cdot)) \in \mathbb R$ is continuous. If there exists a vector $l_0 \in \mathbb R^n$ such that $\partial^-_{l_0}\varphi(t,z,w(\cdot)) = + \infty$, then $\partial^-_l\varphi(t,z,w(\cdot)) + \infty$ for every $l \in \mathbb R^n$.
\end{lemma}
The lemma follows directly from the inclusion $\varphi \in \Phi$ and definition (\ref{direct_deriv:-}).

\begin{lemma}\label{lem_partial_lk_gk}
Let $\varphi \in \Phi$, $(t,z,w(\cdot)) \in \mathbb G$, $t < \vartheta$ and $\tau_k \in [t,\vartheta]$, $l_k, l_*, g_k, g_* \in \mathbb R^n$, $k \in \mathbb N$. Let
$\tau_k \to t$, $l_k \to l_*$, $g_k \to g_*$ as $k \to \infty$. Then
\begin{equation}\label{lem_partial_lk_gk-0}
\partial^-_{l_*} \varphi(t,z,w(\cdot)) \leq \liminf\limits_{k \to \infty} \frac{\varphi(\tau_k,x^{l_k}(\tau_k),x^{g_k}_{\tau_k}(\cdot)) - \varphi(t,z,w(\cdot))}{\tau_k - t}.
\end{equation}
\end{lemma}
The lemma can be proved similar to Lemma \ref{lem_partial_k}.
%

\begin{lemma}\label{len_positive_partial}
Let $\varphi \in \Phi$ and $(t,z,w(\cdot)) \in \mathbb G$, $t < \vartheta$. Let $L \subset \mathbb R^n$ be a nonempty compact set. Suppose that
\begin{equation}\label{len_positive_partial_0}
\partial^-_l \varphi(t,z,w(\cdot)) > 0,\quad l \in L.
\end{equation}
Then there exist numbers $\epsilon_*, \delta_* > 0$ such that
\begin{equation}\label{len_positive_partial_1}
\varphi(\tau,x^l(\tau),x^g_\tau(\cdot)) - \varphi(t,z,w(\cdot)) > \epsilon_* (\tau - t),\quad \tau \in (t,t+\delta_*],\quad l,g \in L.
\end{equation}
\end{lemma}
\proof In the case when there exists $l_0 \in L$ such that $\partial^-_{l_0} \varphi(t,z,w(\cdot)) < +\infty$, according to Lemma \ref{lem_partial_coun} and compactness of $L$, there exists $\epsilon_* > 0$ such that
\begin{equation}\label{len_positive_partial_2}
\epsilon_* < (1/2)\min\limits_{l \in L} \partial^-_l \varphi(t,z,w(\cdot)).
\end{equation}
In the case when $\partial^-_l \varphi(t,z,w(\cdot)) = +\infty$ for every $l \in L$, we take an arbitrary number $\epsilon_* > 0$. For the sake of a contradiction, suppose that, for each $k \in \mathbb N$, there exists a number $\tau_k \in (t,t + 1/k]$ and vectors $l_k,g_k \in L$ such that
\begin{equation}\label{len_positive_partial_3}
\varphi(\tau_k,x^{l_k}(\tau_k),x^{g_k}_{\tau_k}(\cdot)) - \varphi(t,z,w(\cdot)) \leq \epsilon_*(\tau_k - t).
\end{equation}
Since $L$ is compact, then, without loss of generality, we can suppose that $l_k \to l_*$ and $g_k \to g_*$ as $k \to \infty$. Then, applying Lemma \ref{lem_partial_lk_gk}, we obtain $\partial^-_{l_*} \varphi(t,z,w(\cdot)) \leq \epsilon_*$. It contradicts (\ref{len_positive_partial_2}). \hfill $\square$

\begin{lemma}\label{lem_MVI}
Let $\varphi \in \Phi$ and $(t,z,w(\cdot)) \in \mathbb G$, $t < \vartheta$. Let $L \subset \mathbb R^n$ be the nonempty convex compact set. Suppose that {\em (\ref{len_positive_partial_0})} holds. Then, for every $\delta \in (0,\vartheta - t)$, there exist
\begin{equation}\label{lem_MVI-0}
\begin{array}{c}
(\tau_*,v_*) \in \Omega_\delta = \big\{(\tau,v) \in [t,t+\delta] \times \mathbb R^n \colon \min\limits_{l \in L} \|v - z - l(\tau - t)\| \leq \delta\big\},\\[0.2cm]
g_* \in L,\quad (p_0,p) \in D^-\varphi(\tau_*,v_*,x^{g_*}_{\tau_*}(\cdot))
\end{array}
\end{equation}
such that
\begin{equation}\label{lem_MVI-1}
p_0 + \langle l, p\rangle > 0,\quad l \in L.
\end{equation}
\end{lemma}
\proof By the definition of $\Omega_\delta$, one can take $\alpha > 0$ such that $\|v\| \leq \alpha$, for every $(\tau,v) \in \Omega_{(\vartheta - t)}$.
Due to the inclusion $\varphi \in \Phi$, there exists $\lambda_\varphi = \lambda_\varphi(\alpha)$ such that
\begin{equation}\label{lem_MVI-lambdarho}
|\varphi(\tau,v,r(\cdot)) - \varphi(\tau,v,r'(\cdot))| \leq \lambda_\varphi\|r(\cdot) - r'(\cdot)\|_1
\end{equation}
for any $\tau \in [t,\vartheta]$ and $(v,r(\cdot)),(v,r'(\cdot)) \in P(\alpha)$. According to Lemma \ref{len_positive_partial}, there exist $\epsilon_*,\delta_* > 0$ such that (\ref{len_positive_partial_1}) holds. Let $\lambda_L = \max\{\|l\|\,|\,l \in L\}$. Then, without loss of generality, we can suppose that
\begin{equation}\label{lem_MVI-delta}
\delta \leq \delta_*,\quad \delta < \epsilon_* / (\lambda_\varphi (1 + 2\lambda_L)).
\end{equation}


On the set $\overline{\Omega}_\delta = \Omega_\delta \times L \times [t,t + \delta] \times L$, for each $k \in \mathbb N$, let us define the function
\begin{equation*}
\gamma_k(\tau,v,g,\xi,l) = \varphi(\tau,v,x^g_\tau(\cdot)) + k\|v - z - l (\xi - t)\|^2 + k (\tau - \xi)^2 - \epsilon_* (\xi - t),
\end{equation*}
\begin{equation}\label{lem_MVI-gamma_k}
(\tau,v,g,\xi,l) \in \overline{\Omega}_\delta.
\end{equation}
Using Lemma \ref{lem_lim}, one can show that this function is continuous. The set $\overline{\Omega}_\delta$ is compact. Therefore, there exists a point $(\tau_k,v_k,g_k,\xi_k,l_k) \in \overline{\Omega}_\delta$ such that
\begin{equation}\label{tauk_rk_fk_xik}
\gamma_k(\tau_k,v_k,g_k,\xi_k,l_k) = \min\limits_{(\tau,v,g,\xi,l) \in \overline{\Omega}_\delta} \gamma_k(\tau,v,g,\xi,l).
\end{equation}
Furthermore, without loss of generality, we suppose that $(\tau_k,v_k,g_k,\xi_k,l_k) \to (\overline{\tau},\overline{v},\overline{g},\overline{\xi},\overline{l}) \in \overline{\Omega}_\delta$ as $k \to \infty$. Due to (\ref{tauk_rk_fk_xik}), we have
\begin{equation}\label{varphi_k}
\gamma_k(\tau_k,v_k,g_k,\xi_k,l_k) \leq \gamma_k(t,z,g_k,t,l_k) = \varphi(t,z,w(\cdot)).
\end{equation}
Consequently, we obtain
\begin{equation}\label{line_tau_line_xi}
\overline{\tau} = \overline{\xi},\quad \overline{v} = z + \overline{l} (\overline{\xi} - t).
\end{equation}

Let us show that $\overline{\xi} < t+\delta$. For the sake of a contradiction, suppose that $\overline{\xi} = t + \delta$. Then, applying Lemma \ref{lem_lim} and (\ref{len_positive_partial_1}), (\ref{lem_MVI-delta}), (\ref{lem_MVI-gamma_k}), (\ref{line_tau_line_xi}), we derive
\begin{equation*}
\begin{array}{c}
\liminf\limits_{k \to \infty} \gamma_k(\tau_k,v_k,g_k,\xi_k,l_k) \geq \lim\limits_{k \to \infty}\big(\varphi(\tau_k,v_k,x^{g_k}_{\tau_k}(\cdot)) - \epsilon_* (\xi_k - t) \big) \\[0.3cm]
 = \varphi(t+\delta,z+\overline{l}\delta,x^{\overline{g}}_{t+\delta}(\cdot)) - \epsilon_* \delta > \varphi(t,z,w(\cdot)).
 \end{array}
\end{equation*}
This inequality contradicts (\ref{varphi_k}).

In accordance with $\overline{\xi} < t+\delta$ and (\ref{line_tau_line_xi}), one can take $k \in \mathbb N$ so that
\begin{equation}\label{g_k_q_k}
\xi_k < t + \delta,\quad \tau_k < t+\delta,\quad \|v_k - z - l_k (\xi_k - t)\| \leq \delta / 3,\quad \lambda_L |\tau_k - \xi_k| \leq \delta / 3,
\end{equation}
where the number $\lambda_L$ is chosen above. Put
\begin{equation}\label{lem_MVI-star}
p_0 = - \lambda_\varphi \|v_k - z - g_k (\tau_k - t)\| - 2 k (\tau_k - \xi_k),\quad
p = - 2 k (v_k - z - l_k (\xi_k - t)).
\end{equation}

Let us prove the inclusion $(p_0,p) \in D^- \varphi (\tau_k,v_k,x^{g_k}_{\tau_k}(\cdot))$. Take a vector $l \in \mathbb R^n$ and define a function $y^l(\cdot) \in \Lambda(\tau_k,v_k,x^{g_k}_{\tau_k}(\cdot))$ by $y^l(\tau) = v_k + l (\tau - \tau_k)$, $\tau \in [\tau_k,\vartheta]$.
Let $\delta_l = \min\{t - \tau_k + \delta, \delta/(3\|l - l_k\| + 1)\}$. Then, applying (\ref{g_k_q_k}), for every $\tau \in [\tau_k,\tau_k + \delta_l]$, we get
\begin{equation*}
\|y^l(\tau) - z - l_k (\tau - t)\| \leq \|l - l_k\| (\tau - \tau_k) + \|v_k - z - l_k (\xi_k - t)\| + |\xi_k - \tau_k| \|l_k\| \leq \delta.
\end{equation*}
Therefore, taking into account the definition of $\Omega_\delta$ in (\ref{lem_MVI-0}), we have $(\tau,y^l(\tau)) \in \Omega_\delta$, for every $\tau \in [\tau_k,\tau_k + \delta_l]$. Then, applying (\ref{lem_MVI-gamma_k}), (\ref{tauk_rk_fk_xik}), we obtain
\begin{equation}\label{lem_MVI-2}
\begin{array}{c}
0 \leq \gamma_k(\tau,y^{l}(\tau),g_k,\xi_k,l_k) -  \gamma_k(\tau_k,v_k,g_k,\xi_k,l_k) \\[0.3cm]
= \varphi(\tau,y^l(\tau),x^{g_k}_\tau(\cdot)) - \varphi(\tau_k,v_k,x^{g_k}_{\tau_k}(\cdot)) + k \|l\|^2 (\tau - \tau_k)^2 \\[0.3cm]
+ 2 k \langle l, v_k - z - l_k (\xi_k - t) \rangle (\tau - \tau_k) + k (\tau - \tau_k)^2 + 2 k (\tau_k - \xi_*) (\tau - \tau_k).
\end{array}
\end{equation}
Furthermore, using (\ref{lem_MVI-lambdarho}), we derive
\begin{equation}\label{lem_MVI-3}
\begin{array}{c}
|\varphi(\tau,y^l(\tau),y^l_\tau(\cdot)) - \varphi(\tau,y^{l}(\tau),x^{g_k}_\tau(\cdot))| \\[0.3cm]
\leq \lambda_\varphi  \|l - g_k\| (\tau - \tau_k)^2 + \lambda_\varphi  \|v_k - z - g_k (\tau_k - t)\| (\tau - \tau_k).
\end{array}
\end{equation}
From (\ref{lem_MVI-star})--(\ref{lem_MVI-3}), taking into account (\ref{direct_deriv:-}), (\ref{D:-}), we obtain the inclusion $(p_0,p) \in D^-\varphi(\tau_k,v_k,x^{g_k}_{\tau_k}(\cdot))$.

Let us prove (\ref{lem_MVI-1}). We first consider the case when $\xi_k > t$. Let $l \in L$. Since $L$ is convex, then we have $l_\nu = l_k + \nu (l - l_k)/(\xi_k - t) \in L$, $\nu \in [0,\xi_k - t]$. Then, according to (\ref{lem_MVI-gamma_k}), (\ref{tauk_rk_fk_xik}), for every $\nu \in [0,\min\{\xi_k - t,t + \delta - \xi_k\}]$, we derive
\begin{equation*}
\begin{array}{c}
0 \leq \gamma_k(\tau_k,v_k,g_k,\xi_k + \nu,l_\nu) - \gamma_k(\tau_k,v_k,g_k,\xi_k,l_k)\\[0.3cm]
= k \nu^2 \|l + l_\nu - l_k\|^2  - 2 k \nu \langle l + l_\nu - l_k, v_k - z - l_k (\xi_k - t)\rangle\\[0.3cm]
 + k \nu^2 - 2 k (\tau_k - \xi_k) \nu - \epsilon_* \nu.
\end{array}
\end{equation*}
Dividing this inequality by $\nu$ and passing to the limit as $\nu \to 0$, we get
\begin{equation}\label{lem_MVI-4}
\epsilon_* \leq - 2 k\langle l, v_k - z - l_k (\xi_k - t) \rangle - 2 k (\tau_k - \xi_k).
\end{equation}
Using (\ref{lem_MVI-delta}), (\ref{g_k_q_k}), we have
\begin{equation*}
\|v_k - z - g_k (\tau_k - t)\| \leq \|v_k - z - l_k (\xi_k - t)\| + \|l_k\| (\xi_k - t) + \|g_k\| (\tau_k - t)
\end{equation*}
\begin{equation}\label{lem_MVI-5}
\leq (1 + 2\lambda_L) \delta < \epsilon_* / \lambda_\varphi.
\end{equation}
From (\ref{lem_MVI-star}), (\ref{lem_MVI-4}), (\ref{lem_MVI-5}) we obtain (\ref{lem_MVI-1}).

Let us consider the case when $\xi_k = t$. In accordance with (\ref{lem_MVI-gamma_k}), (\ref{tauk_rk_fk_xik}), for every $l \in L$ and $\nu \in [0,\delta]$, we have
\begin{equation*}
\begin{array}{c}
0 \leq \gamma_k(\tau_k,v_k,g_k,t + \nu,l) - \gamma_k(\tau_k,v_k,g_k,t,l_k) \\[0.3cm]
= k \|l\|^2 \nu^2 - 2 k \nu \langle l, v_k - z\rangle + k \nu^2 - 2 k (\tau_k - t) \nu - \epsilon_* \nu.
\end{array}
\end{equation*}
Dividing this inequality by $\nu$ and passing to the limit as $\nu \to 0$, we get
\begin{equation*}
\epsilon_* \leq - 2 k\langle l, v_k - z\rangle - 2 k (\tau_k - t).
\end{equation*}
From this estimate, taking into account (\ref{lem_MVI-star}), (\ref{lem_MVI-5}) we conclude (\ref{lem_MVI-1}). \hfill $\square$

\proof $(d) \Rightarrow (c)$. Let $\varphi$ be a viscosity solution of problem (\ref{HJ}), (\ref{HJ0}). Then $\varphi$ satisfies the inclusion $\varphi \in \Phi$ and terminal condition (\ref{HJ0}).

Let us prove (\ref{ineq_for_deriv:-}). For the sake of a contradiction, suppose that there exist $(t,z,w(\cdot)) \in \mathbb G$, $t < \vartheta$ and $s \in \mathbb R^n$ such that
\begin{equation*}
\partial^-_l \varphi(t,z,w(\cdot)) + H(t,z,w(-h),s) - \langle l,s \rangle > 0,\quad l \in F(z,w(-h)).
\end{equation*}
If there exists $l_0 \in F(z,w(-h))$ such that $\partial^-_{l_0} \varphi(t,z,w(\cdot)) \in \mathbb R$, then, taking into account condition $(H_1)$ and Lemma \ref{lem_partial_coun}, one can take $\eta, \epsilon > 0$ so that
\begin{equation}\label{4to3_1}
\partial^-_l \varphi(t,z,w(\cdot)) + H(t,z,w(-h),s) - \langle l,s \rangle > \epsilon,\quad l \in [F(z,w(-h))]^\eta,
\end{equation}
where the symbol $[F]^\eta$ denote the closed $\eta$-neighborhood of a set $F \subset \mathbb R^n$. If there exists $l_0 \in F(z,w(-h))$ such that $\partial^-_{l_0} \varphi(t,z,w(\cdot)) = + \infty$, then, in accordance with Lemma \ref{lem_partial_coun}, inequality (\ref{4to3_1}) also holds.

Put $L = [F(z,w(-h))]^\eta$. Since $w(\cdot) \in \mathrm{PC}$, then there exists a number $\tau_0 \in (t,\min\{\vartheta,t+h\})$ such that $w(\cdot)$ is continuous on $[-h,\tau_0-t-h]$. Then, according to $(H_1)$, $(F)$, there exists a number $\delta \in (0, \tau_0 - t)$ such that
\begin{equation}\label{4to3_2}
\begin{array}{c}
|H(\tau,v,w(\tau-t-h),s) - H(t,z,w(-h),s)| \leq \epsilon, \\[0.2cm]
F(v,w(\tau-t-h)) \subset [F(z,w(-h))]^\eta,
\end{array}
\quad (\tau,v) \in \Omega_\delta,
\end{equation}
where $\Omega_\delta$ is defined by (\ref{lem_MVI-0}). Define the functional $\tilde{\varphi} \colon \mathbb G \mapsto \mathbb R$ by
\begin{equation*}
\tilde{\varphi}(\tau,v,r(\cdot)) = \varphi(\tau,v,r(\cdot)) + (H(t,z,w(-h),s) - \epsilon)(\tau - t) - \langle v,s \rangle ,\  (\tau,v,r(\cdot)) \in \mathbb G.
\end{equation*}
Science $\varphi \in \Phi$, then the inclusion $\tilde{\varphi} \in \Phi$ holds. Furthermore, from (\ref{4to3_1}), we have $\partial^-_l \tilde{\varphi}(t,z,w(\cdot)) > 0$, $l \in L$. Applying Lemma \ref{lem_MVI} to the functional $\tilde{\varphi}$ and the set $L$, we obtain that there exist $(\tau_*,v_*) \in \Omega_\delta$, $g_* \in L$ and $(p_0,p) \in D^-\tilde{\varphi}(\tau_*,v_*,x^{g_*}_{\tau_*}(\cdot))$ such that
\begin{equation}\label{4to3_3}
p_0 + \langle l,p \rangle > 0,\quad l \in L = [F(z,w(-h))]^\eta.
\end{equation}
Let us define $p'_0 = p_0 - H(t,z,w(\cdot),s) + \epsilon$ and $p' = p + s$. Then we have $(p'_0,p') \in D^-\varphi(\tau_*,v_*,x^{g_*}_{\tau_*}(\cdot))$. Thus, from (\ref{visc_sol:-}), (\ref{4to3_2}) and $(H_3)$, we obtain
\begin{equation*}
0 \geq p'_0 + H(\tau_*,v_*,w(\tau_*-t-h),p') \geq p_0 + \min\limits_{l \in [F(z,w(-h))]^\eta} \langle l,p \rangle,
\end{equation*}
that contradict (\ref{4to3_3}). Thus, (\ref{ineq_for_deriv:-}) holds. For $\partial^+_l \varphi(t,z,w(\cdot))$ and $D^+ \varphi(t,z,w(\cdot))$ one can establish statements similar to Lemmas \ref{lem_partial_coun}--\ref{lem_MVI} and prove (\ref{ineq_for_deriv:+}).\hfill $\square$

%
%

\subsection{Proof $(c) \Rightarrow (b)$}

Let $(t,z,w(\cdot)) \in \mathbb G$ and $\eta \geq 0$. Denote by $\overline{X}(t,z,w(\cdot),\eta)$ the set of the functions $x(\cdot) \in \Lambda(t,z,w(\cdot))$ that satisfy the following delay differential inclusion:
\begin{equation}\label{inclution_eta}
\dot{x}(\tau) \in [F(x(\tau),x(\tau-h))]^\eta \text{ for a.e. } \tau \in [t,\vartheta],
\end{equation}
where the symbol $[F]^\eta$ denote the closed $\eta$-neighborhood of a set $F \subset \mathbb R^n$.

\begin{lemma}\label{lem_X_eta}
Let $(t,z,w(\cdot)) \in \mathbb G$, $t < \vartheta$ and $\eta_*,\eta_k \geq 0$, $k \in \mathbb N$. Let $\eta_k \to \eta_*$ as $k \to \infty$. Let a sequence
$x^{(k)}(\cdot) \in \overline{X}(t,z,w(\cdot),\eta_k)$, $k \in \mathbb N$ be chosen. Then there exists a subsequence $x^{(k_i)}(\cdot)$ and a function $x^*(\cdot) \in \overline{X}(t,z,w(\cdot),\eta_*)$ such that
\begin{equation}\label{lem_X_eta_0}
\max\limits_{\tau \in [t-h,\vartheta]} \|x^{(k_i)}(\tau) - x^*(\tau)\| \to 0 \text{ as } i \to \infty.
\end{equation}
\end{lemma}
\proof By analogy with Proposition \ref{ass_X}, one can show the existence of numbers $\overline{\alpha}_X,\overline{\lambda}_X > 0$ such that, for every $k \in \mathbb N$ and $\tau,\tau' \in [t,\vartheta]$, we have
\begin{equation}\label{lem_X_eta_1}
(x^{(k)}(\tau),x^{(k)}_\tau(\cdot)) \in P(\overline{\alpha}_X),\quad
\|x^{(k)}(\tau) - x^{(k)}(\tau')\| \leq \overline{\lambda}_X|\tau - \tau'|.
\end{equation}
Then, due to C. Arzela--G. Ascoli Theorem (see, e.g., \cite[p. 207]{Natanson_1960}), there exists a subsequence $x^{(k_i)}(\cdot)$ and a function $x^*(\cdot) \in \Lambda(t,z,w(\cdot))$ such that (\ref{lem_X_eta_0}) holds.

The inclusion $x^*(\cdot) \in \overline{X}(t,z,w(\cdot),\eta_*)$ can be proved similar to Lemma 1 of \cite[p. 76]{Filippov_1988}).
\hfill $\square$

\begin{lemma}\label{lem_eta}
For a functional $\varphi \in \Phi$, inequalities {\em (\ref{minmax_sol:v})}, {\em (\ref{minmax_sol:n})} are equivalent the following inequalities:
\begin{subequations}\label{minmax_eta}
\begin{align}
\inf\limits_{x(\cdot) \in \overline{X}(t,z,w(\cdot),\eta)} \omega(t,x(\cdot),\tau,s) \leq \zeta, \label{minmax_eta:v} \\[0.0cm]
\sup\limits_{x(\cdot) \in \overline{X}(t,z,w(\cdot),\eta)} \omega(t,x(\cdot),\tau,s) \geq \zeta \label{minmax_eta:n}
\end{align}
\end{subequations}
for any $(t,z,w(\cdot)) \in \mathbb G$, $t < \vartheta$, $\tau \in (t,\vartheta]$, $s \in \mathbb R^n$ and $\eta,\zeta > 0$, where
\begin{equation}
\begin{array}{rcl}
\omega(t,x(\cdot),\tau,s) & = & \varphi(\tau,x(\tau),x_\tau(\cdot)) - \varphi(t,x(t),x_t(\cdot)) \\[0.3cm]
& + & \displaystyle\int\limits_t^\tau \big(H(\xi,x(\xi),x(\xi-h),s) - \langle \dot{x}(\xi),s \rangle\big)\mathrm{d}\xi.
\end{array}
\end{equation}
\end{lemma}
\proof It is easy to show that (\ref{minmax_eta:v}), (\ref{minmax_eta:n}) follow from (\ref{minmax_sol:v}), (\ref{minmax_sol:n}).

Let us prove that (\ref{minmax_sol:v}) follow from (\ref{minmax_eta:v}). Let $(t,z,w(\cdot)) \in \mathbb G$, $t < \vartheta$, $\tau \in (t,\vartheta]$ and $s \in \mathbb R^n$. Due to (\ref{minmax_eta:v}), for every $k \in \mathbb N$ there exists a function $x^{(k)}(\cdot) \in \overline{X}(t,z,w(\cdot),1/k)$ such that
\begin{equation}\label{1to2_stable}
\omega(t,x^{(k)}(\cdot),\tau,s) \leq 1/k.
\end{equation}
According to Lemma \ref{lem_X_eta}, there exist a subsequence $x^{(k_i)}(\cdot)$, $i \in \mathbb N$ and a function $x^* (\cdot) \in X(t,z,w(\cdot))$ such that (\ref{lem_X_eta_0}) holds. Then, due to Lemma \ref{lem_lim} and condition $(H_1)$, passing to the limit in (\ref{1to2_stable}) as $i \to \infty$, we conclude (\ref{minmax_sol:v}). In the similar way, one can prove that (\ref{minmax_sol:n}) follows from (\ref{minmax_eta:n}).\hfill $\square$

\proof $(c) \Rightarrow (b)$. Let a functional $\varphi \in \Phi$ satisfies (\ref{HJ0}), (\ref{ineq_for_deriv:-}), (\ref{ineq_for_deriv:+}). According to the definition of a minimax solution of problem (\ref{HJ}), (\ref{HJ0}) and Lemma \ref{lem_eta}, for proving of statement $(b)$, it is sufficient to show that $\varphi$ satisfies (\ref{minmax_eta:v}), (\ref{minmax_eta:n}).

For the sake of a contradiction, suppose that there exist $(t,z,w(\cdot)) \in \mathbb G$, $t < \vartheta$, $\overline{\tau} \in (t,\vartheta]$, $s \in \mathbb R^n$ and $\eta, \zeta > 0$ such that
\begin{equation}\label{1to2_1}
\omega(t,x(\cdot),\overline{\tau},s) > \zeta,\quad x(\cdot) \in \overline{X}(t,z,w(\cdot),\eta).
\end{equation}
Define
\begin{equation}\label{1to2_taus}
\begin{array}{c}
\beta(\tau) = \zeta (\tau - t) / (\overline{\tau} - t),\quad \tau \in [t,\vartheta],\\[0.3cm]
\tau_* = \max\Big\{\tau \in [t,\overline{\tau})\,\Big|\,\min\limits_{x(\cdot) \in \overline{X}(t,z,w(\cdot),\eta)} \omega(t,x(\cdot),\tau,s) \leq \beta(\tau)\Big\}.
\end{array}
\end{equation}
The minimum and maximum in this relation are achieved by virtue of Lemmas \ref{lem_lim}, \ref{lem_X_eta}, condition $(H_1)$ and inequality (\ref{1to2_1}). Therefore, there exists a function $x^*(\cdot) \in \overline{X}(t,z,w(\cdot),\eta)$ such that
\begin{equation}\label{1to2_xcirc}
\omega(t,x^*(\cdot),\tau_*,s) \leq \beta(\tau_*).
\end{equation}
For $l \in \mathbb R^n$, let us define the function $y^l(\cdot) \in \Lambda(\tau_*,x^*(\tau_*),x^*_{\tau_*}(\cdot))$ by the rule $y^l(\tau) = x^*(\tau_*) + l (\tau - \tau_*)$, $\tau \in [\tau_*,\vartheta]$. Since $x^*(\cdot) \in \mathrm{PC}([t-h,\vartheta],\mathbb R^n)$, then there exists a number $\delta_0 \in (0,\min\{h,\overline{\tau} - \tau_*\})$ such that the function $x^*(\cdot)$ is continuous on $[\tau_*-h,\tau_*-h + \delta_0]$. Then, in accordance with $(H_1)$, $(F)$, there exists a number $\delta_* \in (0, \delta_0)$ such that
\begin{equation}\label{1to2_coun}
\begin{array}{c}
|H(\tau_*,x^*(\tau_*),x^*(\tau_*-h),s) - H(\xi,y^l(\xi),x^*(\xi-h),s)| \leq \zeta/(2(\overline{\tau} - t)),\\[0.3cm]
F(x^*(\tau_*),x^*(\tau_*-h)) \subset [F(y^l(\xi),x^*(\xi-h))]^\eta.
\end{array}
\end{equation}
for any $\xi \in [\tau_*,\tau_* + \delta_*]$ and $l \in F(x^*(\tau_*),x^*(\tau_*-h))$. Due to (\ref{ineq_for_deriv:-}), there exist a vector $l' \in F(x^*(\tau_*),x^*(\tau_*-h))$ and a number $\tau' \in (\tau_*,\tau_* + \delta_*]$ such that
\begin{equation*}
\begin{array}{c}
\displaystyle\frac{\varphi(\tau',y^{l'}(\tau'),y^{l'}_{\tau'}(\cdot)) - \varphi(\tau_*,x^*(\tau_*),x^*_{\tau_*}(\cdot))}{\tau' - \tau_*} \\[0.3cm]
+ H(\tau_*,x^*(\tau_*),x^*(\tau_*-h),s) - \langle l',s \rangle \leq \zeta/(2(\overline{\tau} - t)).
\end{array}
\end{equation*}
From this estimate, using the first inequality in (\ref{1to2_coun}), we derive
\begin{equation}\label{1to2_y_l}
\omega(\tau_*,y^{l'}(\cdot),\tau',s) \leq \zeta (\tau' - \tau_*)/(\overline{\tau} - t).
\end{equation}
In accordance with the second inequality in (\ref{1to2_coun}), one can define a function $x'(\cdot) \in \overline{X}(t,z,w(\cdot),\eta)$ so that
$x'(\tau) = x^*(\tau)$, $\tau \in [t,\tau_*]$ and $x'(\tau) = y^{l'}(\tau)$, $\tau \in [\tau_*,\tau']$.
Then, from (\ref{1to2_xcirc}), (\ref{1to2_y_l}), taking into account the definition of $\beta(\cdot)$ in (\ref{1to2_taus}), we obtain
$\omega(t,x'(\cdot),\tau',s) \leq \beta(\tau')$, that contradicts the choice of $\tau_*$ in (\ref{1to2_taus}). Thus, (\ref{minmax_eta:v}) has been proved. In the similar way, one can prove (\ref{minmax_eta:n}). \hfill $\square$

\subsection{Proof $(b) \Rightarrow (a)$}

\begin{lemma}\label{lem_lip}
Let $\varphi \in \Phi$ be ci-differentiable at every point $(t,z,w(\cdot)) \in \mathbb G$, $t < \vartheta$. Let $t \in [t_0,\vartheta]$ and $p(\cdot) \in \mathrm{Lip}([t-h,\vartheta],\mathbb R^n)$. Then the function $\omega(\tau) = \varphi(\tau,p(\tau),p_\tau(\cdot))$, $\tau \in [t,\vartheta]$ is Lipschitz continuous and
\begin{equation}\label{lem_lip-0}
\dot{\omega}(\tau) = \partial^{ci}_{t,w} \varphi(\tau,p(\tau),p_\tau(\cdot)) + \langle \dot{p}(\tau), \nabla_z \varphi(\tau,p(\tau),p_\tau(\cdot)) \rangle
\end{equation}
for almost every $\tau \in [t,\vartheta]$.
\end{lemma}
\proof Lipschitz continuity of the function $\omega(\cdot)$ follows directly from the inclusions $\varphi \in \Phi$ and $p(\cdot) \in \mathrm{Lip}([t-h,\vartheta],\mathbb R^n)$. Let $\tau \in (t,\vartheta)$ be a point such that the derivatives $\dot{\omega}(\tau)$ and $\dot{p}(\tau)$ exist. Then, using ci-differentiability of $\varphi$ at the point $(\tau,p(\tau),p_\tau(\cdot))$, we obtain (\ref{lem_lip-0}) at this point.\hfill $\square$

%

For $\lambda > 1$, $\varepsilon_*(\lambda) = e^{-2\lambda(\vartheta- t_0)}$ and $\varepsilon \in (0,\varepsilon_*(\lambda))$, let us define the functional
\begin{equation}\label{mu}
\begin{array}{c}
\mu_\varepsilon^\lambda(t,z,w(\cdot)) = \nu_\varepsilon^\lambda(t) \eta^\lambda_\varepsilon(z,w(\cdot)),\quad (t,z,w(\cdot)) \in \mathbb G,\\[0.3cm]
\nu_\varepsilon^\lambda(t) = (e^{-2\lambda(t - t_0)} - \varepsilon) / \varepsilon,\quad \eta^\lambda_\varepsilon(z,w(\cdot)) = \sqrt{\varepsilon^4 + \|z\|^2} + \lambda \|w(\cdot)\|_1.
\end{array}
\end{equation}
Then the inclusion $\mu_\varepsilon^\lambda \in \Phi$ holds. Furthermore, $\mu_\varepsilon^\lambda$ is ci-differentiable at every point $(t,z,w(\cdot)) \in \mathbb G$, $t < \vartheta$ and
\begin{eqnarray}
& \displaystyle\partial^{ci}_{t,w} \mu_\varepsilon^\lambda(t,z,w(\cdot)) & = - 2 \lambda (\nu^\lambda_\varepsilon(t)+1) \eta^\lambda_\varepsilon(z,w(\cdot)) + \nu^\lambda_\varepsilon(t) \big(\|z\| - \|w(-h)\|\big),\nonumber\\[0.2cm]
& \displaystyle\nabla_z \mu_\varepsilon^\lambda(t,z,w(\cdot)) & = \big(\nu_\varepsilon(t) / \sqrt{\varepsilon^4 + \|z\|^2}\big) z.\label{nu-deriv}
\end{eqnarray}

\begin{lemma}\label{lem_mu}
Let $(t,z,w(\cdot)) \in \mathbb G$. There exists a number $\lambda = \lambda(t,z,w(\cdot)) > 1$ such that, for every $\varepsilon \in (0,\varepsilon_*(\lambda))$ and every $x(\cdot), y(\cdot) \in X(t,z,w(\cdot))$, the functions $p(\tau) = x(\tau) - y(\tau)$, $\tau \in [t-h,\vartheta]$ and $\omega(\tau) = \mu_\varepsilon^\lambda(\tau,p(\tau),p_\tau(\cdot))$, $\tau \in [t,\vartheta]$ are Lipschitz continuous and
\begin{equation}\label{lem_mu-0}
\dot{\omega}(\tau) + |H(\tau,x(\tau),\kappa(\tau),s(\tau)) - H(\tau,y(\tau),\kappa(\tau),s(\tau)) - \langle \dot{p}(\tau),s(\tau) \rangle| \leq 0
\end{equation}
for almost every $\tau \in [t,\min\{t+h,\vartheta\}]$, where
\begin{equation}\label{lem_mu-1}
\kappa(\tau) = w(\tau - t - h),\quad s(\tau) = \nabla_z \mu_\varepsilon^\lambda (\tau,p(\tau),p_\tau(\cdot)).
\end{equation}
\end{lemma}
\proof Due to Proposition \ref{ass_X} and condition $(H_2)$, there is $\lambda_H > 0$ such that
\begin{equation*}
|H(\tau,x(\tau),\kappa(\tau),s) - H(\tau,y(\tau),\kappa(\tau),s)| \leq \lambda_H \|x(\tau) - y(\tau)\| (1 + \|s\|)
\end{equation*}
for any $\tau \in [t,\min\{t+h,\vartheta\}]$, $x(\cdot), y(\cdot) \in X(t,z,w(\cdot))$ and $s \in \mathbb R^n$. Put $\lambda = \lambda_H + 1$. Let $\varepsilon \in (0,\varepsilon_*(\lambda))$ and $x(\cdot),y(\cdot) \in X(t,z,w(\cdot))$. Then, applying Lemma \ref{lem_lip} to the functional $\mu_\varepsilon^\lambda$ and using (\ref{nu-deriv}), we obtain (\ref{lem_mu-0}).\hfill $\square$

\begin{lemma}\label{lem_psi}
Let $\varphi \in \Phi$ and $(t,z,w(\cdot)) \in \mathbb G$. For every $\lambda > 1$ and $\zeta > 0$, there exists a number $\varepsilon = \varepsilon(t,z,w(\cdot),\lambda,\zeta) > 0$ such that the functionals
\begin{subequations}\label{lem_psi-0}
\begin{align}
\begin{array}{rl}\label{lem_psi-0:-}
\psi_-(\tau,v,r(\cdot)) = \min\limits_{y(\cdot) \in X(t,z,w(\cdot))} & \big(\varphi(\tau,y(\tau),y_{\tau}(\cdot)) \\[0.1cm]
& \quad + \mu_\varepsilon^\lambda(\tau,v - y(\tau),r(\cdot) - y_\tau(\cdot))\big), \end{array} \\[0.2cm]
\begin{array}{rl}\label{lem_psi-0:+}
\psi_+(\tau,v,r(\cdot))= \max\limits_{y(\cdot) \in X(t,z,w(\cdot))} & \big(\varphi(\tau,y(\tau),y_{\tau}(\cdot)) \\[0.1cm]
& \quad - \mu_\varepsilon^\lambda(\tau,v - y(\tau),r(\cdot) - y_\tau(\cdot))\big), \end{array}
\end{align}
\end{subequations}
where $(\tau,v,r(\cdot)) \in \mathbb G$, satisfy the inequalities
\begin{subequations}\label{lem_psi-1}
\begin{align}
|\varphi(\tau,x(\tau),x_\tau(\cdot)) - \psi_-(\tau,x(\tau),x_\tau(\cdot))| \leq \zeta \label{lem_psi-1:-},\\[0.1cm]
|\varphi(\tau,x(\tau),x_\tau(\cdot)) - \psi_+(\tau,x(\tau),x_\tau(\cdot))| \leq \zeta \label{lem_psi-1:+}
\end{align}
\end{subequations}
for any $\tau \in [t,\vartheta]$, $x(\cdot) \in X(t,z,w(\cdot))$.
\end{lemma}
\proof Let $\varphi \in \Phi$ and $(t,z,w(\cdot)) \in \mathbb G$.  Take a number $\alpha > 0$ so that $(z,w(\cdot)) \in P(\alpha)$. According to Proposition~\ref{ass_X}, there exists a number $\alpha_X = \alpha_X(\alpha) > \alpha$ such that (\ref{ass_X-0}) holds. Then, due to the inclusion $\varphi \in \Phi$, there exists a number $\lambda_\varphi = \lambda_\varphi(\alpha_X) > 0$ such that
\begin{equation}\label{lem_psi-lambda_Phi}
\begin{array}{rcl}
|\varphi(\tau,x(\tau),x_{\tau}(\cdot)) & - & \varphi(\tau,y(\tau),y_{\tau}(\cdot))|\\[0.2cm]
& \leq & \lambda_\varphi \big(\|x(\tau) - y(\tau)\| + \|x_\tau(\cdot) - y_\tau(\cdot)\|_1\big)
\end{array}
\end{equation}
for any $\tau \in [t,\vartheta]$ and $x(\cdot),y(\cdot) \in X(t,z,w(\cdot))$.

Let $\lambda > 1$ and $\zeta > 0$. Let us choose $\varepsilon_0 \in (0,\varepsilon_*(\lambda))$ and define $\varepsilon > 0$ so that
\begin{equation}\label{lem_psi-epsilon}
\varepsilon < \zeta,\quad \varepsilon < \varepsilon_0,\quad
\varepsilon < \zeta (\varepsilon_*(\lambda) - \varepsilon_0)/ (\lambda_\varphi \theta),\quad \theta = 2 \lambda_\varphi (1 + h) \alpha_X + \zeta.
\end{equation}
We will show that (\ref{lem_psi-1:-}) holds. Let $\tau \in [t,\vartheta]$ and $x(\cdot) \in X(t,z,w(\cdot))$. In accordance with (\ref{mu}), (\ref{lem_psi-0:-}), (\ref{lem_psi-epsilon}), we have
\begin{equation}\label{lem_psi-2}
\psi_-(\tau,x(\tau), x_\tau(\cdot)) - \varphi(\tau,x(\tau),x_\tau(\cdot)) \leq  \mu_\varepsilon^\lambda(\tau,0,r(\cdot) \equiv 0) \leq \zeta.
\end{equation}
Due to Lemma \ref{lem_lim} and Lemma \ref{lem_eta} (for $\eta_k = \eta = 0$), the minimum in (\ref{lem_psi-0:-}) is attained. Therefore, there exists a function $y^*(\cdot) \in X(t,z,w(\cdot))$ such that
\begin{equation}\label{lem_psi-3}
\psi_-(\tau,x(\tau), x_\tau(\cdot)) = \varphi(\tau,y^*(\tau),y^*_\tau(\cdot)) + \mu_\varepsilon^\lambda(\tau,x(\tau)-y^*(\tau), x_\tau(\cdot)-y^*_\tau(\cdot)).
\end{equation}
Hence, taking into account (\ref{ass_X-0}), (\ref{lem_psi-lambda_Phi})--(\ref{lem_psi-2}), we derive
\begin{equation*}
\mu_\varepsilon^\lambda(\tau,x(\tau)-y^*(\tau), x_\tau(\cdot)-y^*_\tau(\cdot)) \leq \varphi(\tau,x(\tau),x_\tau(\cdot)) - \varphi(\tau,y^*(\tau),y^*_{\tau}(\cdot)) + \zeta \leq \theta.
\end{equation*}
From this inequality, in accordance with (\ref{mu}), (\ref{lem_psi-epsilon}), we get
\begin{equation*}
\|x(\tau) - y^*(\tau)\| + \|x_\tau(\cdot) - y^*_\tau(\cdot)\|_1 \leq \eta_\varepsilon^\lambda(x(\tau) - y^*(\tau),x_\tau(\cdot) - y^*_\tau(\cdot)) \leq \zeta/\lambda_\varphi.
\end{equation*}
Therefore, in accordance with (\ref{lem_psi-lambda_Phi}), (\ref{lem_psi-3}), we obtain
\begin{equation}\label{lem_psi-4}
\varphi(\tau,x(\tau),x_\tau(\cdot)) - \psi_-(\tau,x(\tau), x_\tau(\cdot)) \leq \zeta.
\end{equation}
From (\ref{lem_psi-2}), (\ref{lem_psi-4}) we conclude (\ref{lem_psi-1:-}). In the similar way, (\ref{lem_psi-1:+}) can be proved.\hfill $\square$

\proof $(b) \Rightarrow (a)$. Let $\varphi \in \Phi$ be a minimax solution of problem (\ref{HJ}), (\ref{HJ0}). It means that $\varphi$ satisfies (\ref{HJ0}) and (\ref{minmax_sol:v}), (\ref{minmax_sol:n}). For proving statement $(a)$, we need to show that $\varphi$ satisfies (\ref{dynamic_programming}), in which, without loss of generality, we can suppose that $\tau < t+h$.

Let $(t,z,w(\cdot)) \in \mathbb G$, $t < \vartheta$, $\overline{\tau} \in (t,\min\{t+h,\vartheta\})$ and $\zeta > 0$. According to Lemmas \ref{lem_mu}, \ref{lem_psi}, let us choose $\lambda = \lambda(t,z,w(\cdot))$, $\varepsilon = \varepsilon(\zeta / 3)$ and define the functionals $\mu_\varepsilon^\lambda$, $\psi_-$  and $\psi_+$. Then, for proving (\ref{dynamic_programming}), it is sufficient to show the following inequalities:
\begin{subequations}
\begin{align}
\begin{array}{rl}\label{2to1-0:-}
\displaystyle\inf\limits_{u(\cdot) \in \mathfrak{U}(t)} \bigg(\psi_-(\overline{\tau},x(\overline{\tau}),x_{\overline{\tau}}(\cdot)) & + \displaystyle\int\limits_t^{\overline{\tau}} f^0(\xi,x(\xi),x(\xi-h),u(\xi))\mathrm{d}\xi \bigg) \\[0.0cm]
& \quad\quad\quad\quad\quad\quad \leq \psi_-(t,z,w(\cdot)) + \zeta/3,
\end{array} \\[0.0cm]
\begin{array}{rl}\label{2to1-0:+}
\displaystyle\inf\limits_{u(\cdot) \in \mathfrak{U}(t)} \bigg(\psi_+(\overline{\tau},x(\overline{\tau}),x_{\overline{\tau}}(\cdot)) & + \displaystyle\int\limits_t^{\overline{\tau}} f^0(\xi,x(\xi),x(\xi-h),u(\xi))\mathrm{d}\xi \bigg) \\[0.0cm]
& \quad\quad\quad\quad\quad\quad \geq \psi_+(t,z,w(\cdot)) - \zeta/3,
\end{array}
\end{align}
\end{subequations}
where $x(\cdot) = x(\cdot\,|\, t,z,w(\cdot),u(\cdot))$ is the motion of system (\ref{system}).

Let $\zeta_* = \zeta/(30(\overline{\tau} - t))$. In accordance with piecewise continuity of $w(\cdot)$ and Proposition \ref{ass_X}, due to conditions $(f_1)$ and $(H_1)$, there exists a partition
\begin{equation*}
t = \tau_1 < \tau_2 < \ldots < \tau_k < \tau_{k+1} = \overline{\tau}.
\end{equation*}
such that for every $x(\cdot),y(\cdot) \in X(t,z,w(\cdot))$ and $u \in \mathbb U$, we have
\begin{equation}\label{2to1-coun}
\begin{array}{l}
\displaystyle\int\limits_{\tau_i}^{\tau_{i+1}} |f^0(\xi,x(\xi),\kappa(\xi),u) - f^0(\tau_i,x(\tau_i),\kappa(\tau_i),u)|\mathrm{d}\xi \leq \zeta_*, \\[0.2cm]
\displaystyle\int\limits_{\tau_i}^{\tau_{i+1}} |\langle f(\xi,x(\xi),\kappa(\xi),u),s(\xi)\rangle - \langle f(\tau_i,x(\tau_i),\kappa(\tau_i),u),s(\tau_i)\rangle|\mathrm{d}\xi \leq \zeta_*, \\[0.2cm]
\displaystyle\int\limits_{\tau_i}^{\tau_{i+1}} |H(\xi,x(\xi),\kappa(\xi),s(\xi)) - H(\tau_i,x(\tau_i),\kappa(\tau_i),s(\tau_i))|\mathrm{d}\xi \leq \zeta_*, \\[0.2cm]
\displaystyle\int\limits_{\tau_i}^{\tau_{i+1}} |H(\xi,y(\xi),\kappa(\xi),s(\xi)) - H(\xi,y(\xi),\kappa(\xi),s(\tau_i))|\mathrm{d}\xi \leq \zeta_*, \\[0.2cm]
\displaystyle\int\limits_{\tau_i}^{\tau_{i+1}} |\langle \dot{y}(\xi),s(\xi) \rangle - \langle \dot{y}(\xi),s(\tau_i) \rangle|\mathrm{d}\xi \leq \zeta_*,
\end{array}
\end{equation}
where $\kappa(\cdot)$ and $s(\cdot)$ are defined by (\ref{lem_mu-1}).

Define $u(\cdot) \in \mathfrak{U}(t)$ and the corresponding motion $x(\cdot) = x(\cdot\,|\,t,z,w(\cdot),u(\cdot))$ of system (\ref{system}) so that
\begin{equation}\label{2to1_ui}
\begin{array}{c}
u(\tau) = u_i \in \argmin\limits_{u \in \mathbb U} \big(\langle f(\tau_i,x(\tau_i),\kappa(\tau_i),u),s_i \rangle + f^0(\tau_i,x(\tau_i),\kappa(\tau_i),u) \big),\\[0.2cm]
\tau \in [\tau_i,\tau_{i+1}),\quad i \in \overline{1,k},
\end{array}
\end{equation}
where
\begin{equation}\label{2to1_si}
s_i = \nabla_z\mu_\varepsilon^\lambda(\tau_i,p^{(i)}(\tau_i), p^{(i)}_{\tau_i}(\cdot)),\quad
p^{(i)}(\tau) = x(\tau) - y^{(i)}(\tau),\quad \tau \in [t-h,\vartheta],
\end{equation}
and the function $y^{(i)}(\cdot) \in X(t,z,w(\cdot))$ satisfies the equality
\begin{equation}\label{2to1_yi}
\varphi(\tau_i,y^{(i)}(\tau_i),y^{(i)}_{\tau_i}(\cdot)) + \mu_\varepsilon^\lambda(\tau_i,p^{(i)}(\tau_i),p^{(i)}_{\tau_i}(\cdot))
= \psi_-(\tau_i,x(\tau_i),x_{\tau_i}(\cdot)).
\end{equation}
The minimum in (\ref{2to1_ui}) is attained according to condition $(f_1)$ and a compactness of the set $\mathbb U$.

Let $i \in \overline{1,k}$. Take the vector $s_i$ from (\ref{2to1_si}). Due to (\ref{minmax_sol:v}), there exists $y^{(i)}(\cdot) \in X(t,z,w(\cdot))$ that satisfies (\ref{2to1_yi}) and the inequality
\begin{equation*}
\begin{array}{rcl}
\varphi(\tau_{i+1},y^{(i)}(\tau_{i+1}),y^{(i)}_{\tau_{i+1}}(\cdot)) & +& \displaystyle\int\limits_{\tau_i}^{\tau_{i+1}} \big(H(\xi,y^{(i)}(\xi),\kappa(\xi),s_i) - \langle \dot{y}^{(i)}(\xi),s_i \rangle\big)\mathrm{d}\xi\\[0.2cm]
& \leq & \varphi(\tau_i,y^{(i)}(\tau_i),y^{(i)}_{\tau_i}(\cdot)) + \zeta (\tau_{i+1} - \tau_i) / (6(\overline{\tau} - t)).
\end{array}
\end{equation*}
From this inequality, applying Lemma \ref{lem_mu} and (\ref{H}), (\ref{lem_psi-0:-}), (\ref{2to1-coun})--(\ref{2to1_yi}), we derive
\begin{equation*}
\begin{array}{rcl}
\psi_-(\tau_{i+1},x(\tau_{i+1}),x_{\tau_{i+1}}(\cdot)) & + & \displaystyle\int\limits_{\tau_i}^{\tau_{i+1}} f^0(\xi,x(\xi),\kappa(\xi),u_i)\mathrm{d}\xi \\[0.2cm]
& \leq&  \psi_-(\tau_i,x(\tau_i),x_{\tau_i}(\cdot)) + \zeta (\tau_{i+1} - \tau_i) / (3(\overline{\tau} - t)),
\end{array}
\end{equation*}
and, consequently, we get (\ref{2to1-0:-}).

Let us prove (\ref{2to1-0:+}). Let $u(\cdot) \in \mathfrak{U}(t)$ and $x(\cdot) = x(\cdot\,|\,t,z,w(\cdot),u(\cdot))$. Let $i \in \overline{1,k}$. Due to (\ref{minmax_sol:n}), (\ref{lem_psi-0:+}), there exists $y^{(i)}(\cdot) \in X(t,x,w(\cdot))$ such that
\begin{equation}\label{2to1-yi}
\varphi(\tau_i,y^{(i)}(\tau_i),y^{(i)}_{\tau_i}(\cdot)) - \mu_\varepsilon^\lambda(\tau_i,p^{(i)}(\tau_i),p^{(i)}_{\tau_i}(\cdot)) = \psi_+(\tau_i,x(\tau_i),x_{\tau_i}(\cdot)),
\end{equation}
\begin{equation*}
\begin{array}{rcl}
\displaystyle\varphi(\tau_{i+1},y^{(i)}(\tau_{i+1}),y^{(i)}_{\tau_{i+1}}(\cdot)) & + & \displaystyle\int\limits_{\tau_i}^{\tau_{i+1}} \big(H(\xi,y^{(i)}(\xi),\kappa(\xi),-s_i) + \langle \dot{y}^{(i)}(\xi),s_i \rangle\big)\mathrm{d}\xi \\[0.2cm]
& \geq & \varphi(\tau_i,y^{(i)}(\tau_i),y^{(i)}_{\tau_i}(\cdot)) - \zeta (\tau_{i+1} - \tau_i) / (6(\overline{\tau} - t)),
\end{array}
\end{equation*}
where $s_i$ and $p^{(i)}(\cdot)$ are defined by (\ref{2to1_si}). Then, from Lemma \ref{lem_mu} and (\ref{H}), (\ref{lem_psi-0:+}), (\ref{2to1-coun})--(\ref{2to1_yi}), we derive
\begin{equation*}
\begin{array}{rcl}
\psi_+(\tau_{i+1},x(\tau_{i+1}),x_{\tau_{i+1}}(\cdot)) & + & \displaystyle\int\limits_{\tau_i}^{\tau_{i+1}} f^0(\xi,x(\xi),\kappa(\xi),u(\xi))\mathrm{d}\xi \\[0.2cm]
& \geq & \psi_+(\tau_i,x(\tau_i),x_{\tau_i}(\cdot)) -  \zeta (\tau_{i+1} - \tau_i) / (3(\overline{\tau} - t)),
\end{array}
\end{equation*}
and, consequently, we get (\ref{2to1-0:+}).\hfill $\square$

\section{Conclusions}

On the space of piecewise continuous functions, an optimal control problem for time-delay systems with discrete delay is considered. For the value functional, the corresponding Hamilton-Jacobi-Bellman equation with coinvariant derivatives is investigated. Definitions of minimax and viscosity solutions of the Cauchy problem for this equation are studied. It is proved that both of these solutions exist, are unique and coincide with the value functional. The proof of the viscosity solution uniqueness is based on an analogue of the theorem about "Mean value inequality" for functionals defined on the space of piecewise continuous functions. 

In the future, we plan to get analogous results for more general dynamical systems which motions are described by functional differential equations of delay and neutral types. We plan to consider positional differential games for such systems and apply minimax and viscosity approaches to the corresponding Hamilton-Jacobi equations with a non-convex Hamiltonian. Also, it seems interesting to investigate generalized solutions of boundary value problems for Hamilton-Jacobi equations that correspond to time-optimal control problems and infinite horizon optimal control problems in time-delay systems. Perhaps, results of this paper will be useful for such investigations.

\begin{acknowledgements}
???
\end{acknowledgements}



\end{document}